\documentclass{article}
\usepackage{color} %Color definition
\usepackage{hyperref} %Internal and external links
\usepackage{graphicx} %Graphics inclusion

\usepackage[english]{babel}
\usepackage[nice]{nicefrac}
\usepackage{enumerate}
\usepackage{amsmath,amssymb,amsthm,mathabx}
\usepackage[top=2cm,bottom=2cm,left=2cm,right=2cm]{geometry}

\newtheorem{theorem}{Theorem}[section]                                          
\newtheorem{proposition}[theorem]{Proposition}                          
\newtheorem{lemma}[theorem]{Lemma}

\newtheorem{definition}[theorem]{Definition}
\newtheorem{remark}[theorem]{Remark}

\def\EgalLoi{{~\mathop{= }\limits^{(law)}}~}

\def\supp{{\rm Supp}}
\newcommand{\RRe}{{\mathrm{Re}}}

\newcommand{\Ker}{\mathtt{Ker}}
\def\summ#1#2{{\displaystyle\sum _{\doubleindice{#1}{#2}}}}
\def\doubleindice#1#2{\genfrac{}{}{0pt}{1}{#1}{#2}}

\newcommand{\bs}[1]{\ensuremath{\boldsymbol{#1}}}
\newcommand{\Supp}{\mathtt{Supp}}

\title{Balanced multicolour P{\'o}lya urns via smoothing systems analysis}

\author{C{\'e}cile Mailler\thanks{Department of Mathematical Sciences,
University of Bath,
Claverton Down,
BA2 7AY Bath
United-Kingdom. Email: c.mailler@bath.ac.uk}}

\begin{document}

\maketitle
\begin{abstract}
The present paper aims at describing in details the asymptotic composition of a class of $d$-colour P{\'o}lya urns: namely balanced, tenable and irreducible urns. We decompose the composition vector of such urns according to the Jordan decomposition of their replacement matrix. The projections of the composition vector onto the so-called \emph{small} Jordan spaces are known to be asymptotically gaussian, but the asymptotic behaviour of the projections onto the \emph{large} Jordan spaces are not known in full details up to now and are described by a limit random variable called $W$, depending on the parameters of the urn.

We prove, via the study of smoothing systems, that the variable $W$ has a density and that it is moment-determined. 
\end{abstract}

%\vspace{-\baselineskip}
\tableofcontents

\section{Introduction}
A P{\'o}lya urn is a discrete time stochastic process which was originally introduced by P{\'o}lya and Eggenberger to model the spread of epidemics~\cite{EP23}. Since then, they have been useful in many different areas of mathematics and theoretical computer science and are therefore broadly studied. We can for example cite applications to the analysis of random trees (AVL\footnote{The AVL is a balanced search tree named after its inventors, Adelson-Velskii \& Landis (see~\cite{AVL}).}~\cite{Mahmoud98}, 2--3 trees~\cite{FGP05}), to the analysis of the Bandit algorithm~\cite{LPT04}, or to the reinforced random walks (see for example the survey of ~\cite{Pemantle07}).

The range of methods used to study this random object is also very large. Historically studied by enumerative combinatorics (see for example~\cite{BP85}), P{\'o}lya urns have been efficiently studied by embedding in continuous time (or \emph{poissonization}) since the works of Athreya and Karlin~\cite{AK68} (see for example~\cite{Janson04}). In parallel, since the seminal paper by Flajolet, Dumas \& Puyhaubert~\cite{FDP06}, the analytic combinatorics community successfully tackles the problem in the two-colour case.

A P{\'o}lya urn process is defined as follows: an urn contains balls of different colours, let us denote by $d$ the number of different colours available. 
Fix an initial composition $\bs \alpha= {}^t \!(\alpha_1, \ldots, \alpha_d)\in\mathbb{N}^{d}$, meaning that there are $\alpha_i$ balls of colour $i$ at time zero in the urn, for all $i\in\{1, \ldots, d\}$.
Fix a $d\times d$ matrix $R = (a_{i,j})_{1\leq i,j \leq d}$ with integer coefficients. 
At each step of the process, pick a ball uniformly at random in the urn, 
denote its colour by $c$ ($c\in\{1, \ldots, d\}$), 
put this ball back in the urn and add into the urn $a_{c,i}$ balls of colour $i$ for all $1\leq i\leq d$.
The standard question asked is ``how many balls of each colour are there in the urn'' at time $n$? when $n$ tends to infinity?

Of course, the answer depends on the initial composition vector and on the replacement matrix chosen, and many different behaviours are exhibited in the literature. The aim of the present paper is to investigate into more detail the asymptotic behaviour of a very large class of P{\'o}lya urns: the $d$-colour, balanced, tenable and irreducible urns. This framework is very general and includes, among others, P{\'o}lya urns modelling m-ary trees~\cite{CH01,CP04}, paged binary trees~\cite{CH01,Mah02} or $B$-trees~\cite{CGPT++}. Our precise assumptions are the following:
\begin{itemize}
\item[$\mathtt{(T)}$] The non-diagonal coefficients of $R$ are non-negative, and, for all $1\leq i\leq d$, either $a_{i,i}\geq -1$ or $-a_{i,i}$ is the gcd of $\{a_{1,i}, \ldots, a_{d,i}, \alpha_i\}$. 
\item[$\mathtt{(B)}$] The urn is balanced, meaning that there exists an integer $S$, called the balance, such that, for all $c\in\{1, \ldots, d\}$, $\sum_{i=1}^d a_{c,i} = S$.
\item[$\mathtt{(I)}$] The replacement matrix $R$ is irreducible, meaning that for all $1\leq c, i\leq d$, there exists $n\geq 0$ such that $(R^n)_{c,i}>0$.
\end{itemize}

These assumptions ensure us that $S$ is the eigenvalue of $R$ having maximal real part.
We then consider the Jordan decomposition of the matrix $R$ and fix one Jordan block: 
this Jordan block is associated to a stable subspace~$E$ (containing a unique eigenline) 
and to an eigenvalue $\lambda$. Note that several Jordan blocks can correspond to the same eigenvalue.
We are interested in the behaviour of the projection of the urn composition vector on $E$, 
asymptotically when $n$ tends to infinity. When $\sigma:=\nicefrac{\RRe \lambda}{S} \leq \nicefrac12$. 
This projection, appropriately renormalised, converges in law to a Gaussian distribution, 
independent of the initial composition $\bs \alpha$ of the urn. 
However in the case of ``large'' eigenspaces, i.e. if $\sigma >\nicefrac12$, 
we observe a different behaviour, which we aim at describing better in the present paper.

The projection on a large eigenspace $E$ of the composition vector of the urn at time $n$, renormalised appropriately, converges almost surely and in all $L^p (p\geq 1)$ to a complex random variable $W^{DT}_{\bs \alpha}$ times an eigenvector $v$ associated to $E$. Moreover, if we embed the urn process in continuous time, the obtained Galton-Watson process $U^{CT}(t)$, projected onto $E$, and correctly renormalised, converges almost surely and in all $L^p (p\geq 1)$ to a complex random variable $W^{CT}_{\bs \alpha}$ times $v$. Our main goal is to gather information about the two random variables $W_{\bs \alpha}^{DT}$ and $W_{\bs \alpha}^{CT}$.

The discrete time and the continuous time process are closely related and one of our main tools will be to transport information from one setting to the other. The three main results of the paper are: for all initial composition $\bs \alpha$,
\begin{itemize}
\item the support of the random variables $W^{DT}_{\bs \alpha}$ and $W^{CT}_{\bs \alpha}$ is the whole complex plane if $\lambda\in\mathbb C\setminus \mathbb R$ and the whole real line if $\lambda\in\mathbb R$;
\item $W^{DT}_{\bs \alpha}$ and $W^{CT}_{\bs \alpha}$ both admit a density on $\mathbb C$ (resp. $\mathbb R$ if $\lambda\in\mathbb R$);
\item $W^{DT}_{\bs \alpha}$ and $W^{CT}_{\bs \alpha}$ are moment-determined (the Laplace transform of $W^{DT}_{\bs \alpha}$ converges on the whole plane).
\end{itemize}

These results are the first results in the literature about the variables $W$ induced by multi-colour urns.

To prove them, we first use the branching property of the urn process to prove that it is enough to consider the $d$ atomic initial compositions $\bs e_c$ ($c\in\{1, \ldots, d\}$) defined as follows: all the coefficients of $\bs e_c$ are zero except the $c$th which is equal to 1 if $a_{c,c} \geq 0$, and to $-a_{c,c}$ otherwise. The branching property of urn schemes is largely used in the literature; the novelty of our approach is that we can include the case of P{\'o}lya urn schemes with diagonal coefficients possibly less than~$-1$.
Using again the branching property, we prove that $(W^{DT}_{\bs e_1}, \ldots, W^{DT}_{\bs e_d})$
(resp. $(W^{CT}_{\bs e_1}, \ldots, W^{CT}_{\bs e_d})$) is solution of a system of smoothing equations. 

To our knowledge, this article contains the first analysis of such a smoothing system: 
the literature contains examples of smoothing systems of two equations with real random variable solutions (see~\cite{CMP13} for an example) 
or one smoothing equation with complex solutions (see e.g.~\cite{CP04}), 
but no smoothing system with complex solutions as in the present article.
Such an analysis needed the development of new arguments 
and this approach would be useful in any other context where such smoothing systems would appear.
In a very recent article (appeared on arxiv.org while the current paper was under review), 
\cite{Leckey} gives a survey about smoothing systems and their applications, as well as general results about them.

From the smoothing system in discrete time, we deduce a system of equations verified by the Fourier transforms of the $W^{DT}$s, and prove that these Fourier transforms are integrable, inducing the existence of densities. We then show how to deduce the same result for the continuous time $W$s, and for any initial composition. 

From the system in continuous time, we deduce an induction formula for the moments of the $W^{CT}$s. We then prove by induction that these moments are small enough to apply Carleman's criterion to conclude that they are moment-determined. We then transport the result in discrete time and for all initial initial composition. 

{\bf Plan of the paper:}
We first describe in Section~\ref{sec:preliminaires} our framework and the state of the art concerning the asymptotic behaviour of the multi-colour urns; at the end of Section~\ref{sec:preliminaires}, we also state our main results (namely Theorems~\ref{thm:dcoul_moments} and~\ref{th:density}). 
Section~\ref{sec:CT_arbo} is devoted to characterise the random variables $W^{CT}$ as a solution of a system of fixed point equations in law. In Section~\ref{sec:dcoul_moments}, we use these systems to study the moments of the $W$'s and prove that they are moment-determined, both in continuous and in discrete time. In Section~\ref{sec:discrete}, we imply, from the continuous time systems and from the moment study, that the random variables $W^{DT}$ are also characterised as a solution of a smoothing system. Finally, Section~\ref{sec:d_coul_density} contains the proof that the $W$'s admit a density, both in discrete and continuous time. It is very interesting to see how we will travel from discrete to continuous time all along the paper and how going from one world to the other is very fruitful. But first, let us discuss our set of hypothesis: $\mathtt{(T)}$,~$\mathtt{(B)}$ and~$\mathtt{(I)}$.

\subsection{Discussion of the hypothesis}

Different hypothesis are made in the literature, in order to control the behaviour of P{\'o}lya urns:
\begin{itemize}
\item[$\mathtt{(T_{-1})}$] The coefficients of $R$ are all non-negative except the diagonal coefficients which can be equal to~$-1$.
\item[$\mathtt{(T)}$] The non-diagonal coefficients of $R$ are non-negative, and, for all $1\leq i\leq d$, either $a_{i,i}\geq -1$ or $-a_{i,i}$ is the gcd of $\{a_{1,i}, \ldots, a_{d,i}, \alpha_i\}$. 
\item[$\mathtt{(B)}$] The urn is balanced, meaning that there exists an integer $S$, called the balance, such that, for all $c\in\{1, \ldots, d\}$, $\sum_{i=1}^d a_{c,i} = S$.
\end{itemize}

Assuming $\mathtt{(T_{-1})}$ or $\mathtt{(T)}$ permits to avoid non-tenable urns, i.e. urn schemes in which something impossible is asked: an example would be that you must subtract $3$ balls of colour 1 from the urn while there is only 1 such ball in the urn. 
Allowing only non-negative coefficients in the replacement matrix and coefficients at least $-1$ on the diagonal permits to ensure that the urn is tenable. It is proven in~\cite{Pouyanne08} that assumption $\mathtt{(T)}$ has the same effect although it is much weaker and allows to include a much wider class of P\'olya urns into our framework.
Some authors prefer working without such an assumption but then make all reasoning conditioned to tenability, 
i.e.\ conditioned on the event ``no impossible configuration occur'' (see for example~\cite[Remark~4.2]{Janson04}).

The balance assumption $\mathtt{(B)}$ is quite standard in the literature, though it is not always necessary (this assumption is not needed in~\cite{Janson04}, for example). This assumption implies that the total number of balls in the urn is a deterministic function of time and this property is the foundation of combinatorics approaches while continuous time analysis of P{\'o}lya urns can be done for non-balanced urns.

Finally, we will assume that the urn is irreducible, meaning that any colour can be produced from an initial composition with one unique ball: for a two--colour urn, being irreducible means having a non-triangular replacement matrix. The following definitions define the notion precisely.
\begin{definition}[{see~\cite[page~4]{Janson04}}]
Let $c, i\in\{1, \ldots, d\}$, we say that $c$ {\bf dominates} $i$ 
if there exists $n\geq 1$ such that $(R^n)_{c,i}>0$.
A colour $c\in\{1, \ldots, d\}$ is said to be {\bf dominating} if it dominates every other colour in $\{1, \ldots, d\}$.
\end{definition}

\begin{definition}[{see~\cite[page~4]{Janson04}}]\label{df:dominating}
We say that an urn of replacement matrix $R$ is {\bf irreducible} if and only if every colour is a dominating colour.
\end{definition}

\begin{itemize}
\item[$\mathtt{(I)}$] The replacement matrix $R$ is irreducible.
\end{itemize}

Note that it is sometimes enough (see, for example,~\cite{Janson04} and~\cite{Pouyanne08}) 
to only assume the following weaker version of the irreducibility assumption.
First note that the domination relation is transitive and reflexive and 
thus partitions the set of colours into some equivalence classes:
\begin{definition}[{see~\cite[page~4]{Janson04}}]\label{df:classes}
We say that two colours~$i$ and~$j$ are in the same class if $i$ dominates $j$ and $j$ dominates $i$.
A {\bf dominating class} is a class of dominating colours.
An eigenvalue~$\lambda$ belongs the one of the equivalence class $D$ 
if the restriction of~$R$ to $D$ admits $\lambda$ as an eigenvalue.
\end{definition}
We are now ready to assume this weaker version of irreducibility, which we call $\mathtt{(S)}$ for \emph{simplicity}:
\begin{itemize}
\item[$\mathtt{(S)}$] 
The largest real eigenvalue $\lambda_{max}$ of~$R$ is positive and is a simple eigenvalue of~$R$. 
Furthermore, there is at least one ball of a dominating colour in the urn at time~0, and $\lambda_{max}$ belong to the dominating class. 
\end{itemize}
Note that under $\mathtt{(B)}$, $S=\lambda_{max}$. We refer the reader to Janson~\cite[page~5]{Janson04} for a discussion of this weaker assumption.

Cases of non-irreducible urns are studied in the literature: the diagonal case $R = S I_d$ is the original P{\'o}lya-Eggenberger process and its behaviour is well described in~\cite{Athreya69, BK64, JK77, CMP13}, and triangular urn schemes are developped in~\cite{Janson06, BDM09}. 

Note that a large two-colour urn cannot have negative diagonal coefficients, whereas there exists $d$--colour P{\'o}lya urns with possibly negative coefficients having large eigenvalues. We thus have to include such urns in our study, that is why we only assume $\mathtt{(T)}$ and not the more restrictive assumption $\mathtt{(T_{-1})}$.
In the present paper, we thus choose to assume $\mathtt{(T)}$, $\mathtt{(B)}$ and $\mathtt{(S)}$ (although it is necessary to assume $\mathtt{(I)}$ for Theorem~\ref{th:density}): we can cite many examples of urn processes that fall in this framework (see for example m-ary trees~\cite{CH01,CP04}, paged binary trees~\cite{Mah02,CH01}, $B$-trees~\cite{CGPT++}) and we will see all along the paper how each of these assumptions is used in the proofs. We are interested in the asymptotic behaviour of an urn under these three conditions.

The present setting is different from~\cite{Janson04}'s setting for $d$--colour P{\'o}lya urns ($d\ge 3$), 
where $\mathtt{(T_{-1})}$ is assumed with further assumptions, but~$\mathtt{(B)}$ is not: 
it is however mentioned in~\cite[Remark~4.2]{Janson04} that Janson's main results hold under our set of hypothesis. 
We will thus be able to apply Janson's results before going further in the study of the asymptotic behaviour of the urn.

The following section is devoted to summarising the results of the literature needed as preliminaries to state our main results.

\section{Preliminaries and statement of the main results}\label{sec:preliminaires}
The behaviour of the urn process is already quite well known: We recall hereby the main results of the literature (mainly by \cite{Athreya69}, \cite{Janson04} and \cite{Pouyanne08}) and  thereafter state our main results.

\subsection{Jordan decomposition}
In view of $\mathtt{(B)}$ and $\mathtt{(S)}$, the matrix $R$ admits $S$ as a simple eigenvalue, and every other eigenvalue $\lambda$ of $R$ verifies $\RRe \lambda < S$. The matrix $R$ can be written on its Jordan normal form, meaning that it is similar to a diagonal of blocks $\mathtt{diag}(J_1,\ldots,J_r)$ where each $J_i$ is a matrix shaped as follows:
\[J = \begin{pmatrix}
\lambda & 1 & 0 & \dots & 0\\ 
0 &\lambda & 1 & \ddots & \vdots \\
\vdots & \ddots & \ddots & \ddots & 0\\
\vdots &\dots &\ddots & \lambda & 1\\
0 & \dots & \dots & 0 & \lambda
\end{pmatrix},\]
where $\lambda$ is an eigenvalue of $R$. 
Note that several Jordan blocks can be associated to the same eigenvalue.
In the following, we chose a Jordan block and study the behaviour of the projection of the composition vector onto the subspace associated to this Jordan block. Note that the fact that $R$ is irreducible implies that $S$ is a single eigenvalue of $R$.

\begin{definition}
Let $\lambda$ be an eigenvalue for $R$ and $\sigma = \nicefrac{\RRe\lambda}{S}$. We call $\lambda$ a {\bf large eigenvalue} of $R$ if $\nicefrac12 < \sigma < 1$, or a {\bf small eigenvalue} if $\sigma\leq \nicefrac12$.

A {\bf large Jordan block} is a Jordan block of $R$ associated with a large eigenvalue of $R$ and a {\bf small Jordan block} is a Jordan block associated with a small eigenvalue of $R$.
\end{definition}
We denote by $U_{\bs\alpha}(n) = U(n)\in\mathbb{N}^d$ the composition vector of the urn at time $n$ 
(the subscript $\bs\alpha$ is the initial composition of the urn): 
its $i^{\text{th}}$ coordinate is by definition equal to the number 
of balls of colour $i$ at time $n$ in the urn.
We are interested in the behaviour of $U(n)$ when $n$ tends to infinity. It is showed in the literature that $U(n)$ is easier to describe when decomposed according to the Jordan block decomposition of $R$. 
For every stable subspace~$E$ associated to a Jordan block of $R$, 
we will denote by $\pi_E$ the projection on~$E$ relative to the direct sum of all Jordan subspaces of $R$, 
and we will study separately each projection on a Jordan subspace $E$.

It is standard to embed urn processes in continuous time (see for example~\cite{AK68}): each ball is seen as a clock that rings after a random time with exponential law of parameter one, independently from other clocks in the urn. When a clock rings, it splits into $a_{i,j}+\delta_{i,j}$ balls of colour $j$ ($\forall j\in\{1,\ldots, d\}$) if the clock had colour $i$ (where we uses Kronecker's notation: $\delta_{i,j}$ equals $1$ if $i=j$ and $0$ otherwise). We denote by $\tau_n$ the time of the $n^{\text{th}}$ ring in the urn and by $U^{CT}(t)$ the composition vector of the (continuous time) urn at time~$t$. We have the following standard connection: almost surely,
\begin{equation}\label{eq:dcoul_connexion}
(U(n))_{n\geq 0} = (U^{CT}(\tau_n))_{n\geq 0}.
\end{equation}
In addition, the process $(U(n))_{n\geq 0}$ is independent of the sequence of stopping times $(\tau_n)_{n\geq 0}$.

The asymptotic behaviour of the different projections of $U(n)$ and $U^{CT}(t)$ is partially described in the literature:
\begin{itemize}
\item In continuous time (see~\cite{Janson04}),
\begin{itemize}
\item small projections have a Gaussian behaviour, and
\item renormalised large projections converge almost surely to a random variable $W^{CT}$.
\end{itemize}
\item In discrete time,
\begin{itemize}
\item if $R$ has only small eigenvalues apart from $S$, if $E$ is one of the largest Jordan block associated to the eigenvalue $\lambda$ realising the second highest real part (after $S$), and if $\sigma = \nicefrac{\RRe \lambda}{S}=\nicefrac12$ then projections onto $E$ have a Gaussian behaviour (see \cite[Theorems~3.22 et~3.23]{Janson04}); and
\item renormalised large projections converge almost surely to a random variable $W^{DT}$ (see \cite{Pouyanne08}).
\end{itemize}
\end{itemize}
As one can see, the behaviour of small projections (i.e.\ projections on a small Jordan block) 
in discrete time is not known yet in full generality: Subsection~\ref{sub:small}
is devoted to describing the behaviour of $\pi_E(U(n))$ for all small Jordan space of $R$.
However, the main aim of the paper is to describe the unexplored $W^{DT}$ and $W^{CT}$: 
Subsection~\ref{sub:large} will state the results concerning the projections on large Jordan spaces 
(i.e.\ associated to a large Jordan block), as a preliminary to our main results.

\subsection{Projections on small Jordan spaces}\label{sub:small}
As explained above, the behaviour of the projections of the composition vector (in discrete time) onto the small Jordan blocks is not known yet in full generality. Although our main aim is to focus on the less understood projections onto large Jordan blocks, we sketch here a proof of a general result for small Jordan blocks in order to complete the theory of P{\'o}lya urns under $\mathtt{(B)}$, $\mathtt{(T)}$ and $\mathtt{(S)}$.

\begin{theorem}\label{th:small_lambda}
Under assumptions $\mathtt{(B)}$, $\mathtt{(T)}$ and $\mathtt{(S)}$,
if $E$ is a block of size $\nu+1$ associated to a small eigenvalue $\lambda$ of $R$, 
then there exists a covariance matrix $\Sigma$ such that
\begin{itemize}
\item If $\RRe \lambda = \frac{S}{2}$ then
\[\frac{\pi_E(U_{\bs \alpha}(n))}{\sqrt{Sn\ln^{2\nu+1} n}} \to \mathcal{N}(0,\Sigma),\]
in distribution, asymptotically when $n$ tends to infinity.
\item If $\RRe \lambda < \frac{S}{2}$ then
\[\frac{\pi_E(U_{\bs \alpha}(n))}{\sqrt{Sn}} \to \mathcal{N}(0,\Sigma),\]
in distribution, asymptotically when $n$ tends to infinity.
\end{itemize}
Moreover, $\Sigma$ does not depend on $\bs \alpha$.
\end{theorem}

Let $E$ and $\lambda$ as in Theorem~\ref{th:small_lambda}.
The following result by Janson is the key of the proof of Theorem~\ref{th:small_lambda}:
\begin{theorem}[{\cite[Theorem~3.15 (i) and (ii)]{Janson04}}]\label{th:Janson_small}
Under assumptions $\mathtt{(B)}$, $\mathtt{(T)}$ and $\mathtt{(S)}$,
for all vector $\bs b\in\mathbb{R}^{d}$, define 
\[\tau_{\bs b}(n)= \inf\{t\geq 0 \;|\; \langle \bs b, U^{CT}_{\bs \alpha}(t)\rangle \geq n\}.\]
Then,
\begin{enumerate}[(i)]
\item If $\RRe\lambda = \frac{S}{2}$, then
\[\frac{1}{\sqrt{Sn \ln^{2\nu+1} n}}\;\pi_E(U_{\bs \alpha}^{CT}(\tau_{\bs b}(n))) \to \mathcal N(0,\sigma),\]
in distribution, where $\sigma$ is a covariance matrix.
\item If $\RRe\lambda < \frac{S}{2}$, then
\[\frac1{\sqrt{Sn}}\;\pi_E(U_{\bs \alpha}^{CT}(\tau_{\bs b}(n))) \to \mathcal N(0,\sigma),\]
n distribution, where $\sigma$ is a covariance matrix.
\end{enumerate}
\end{theorem}

Theorem 3.15 in~\cite{Janson04} is slightly different than the above version. The above version corresponds to the special case $z=n$ in Janson's Theorem~3.15. The $(i)$ of \cite[Theorem~3.15]{Janson04} concerns the projection on the union of the small Jordan spaces, and it implies the $(i)$ above by projection on a specified small Jordan space. The $(ii)$ in Theorem~3.15~\cite{Janson04} is more general than the above version, which is the special case $k=\nu$ of Janson's result. The matrix $\sigma$ is given by Equations~(3.11) and (2.15) in~\cite{Janson04} for case $(i)$ above, and by Equations~(3.12) and~(2.16) in~\cite{Janson04} for case $(ii)$ above. It is important to note that $\sigma$ does not depend on $\bs \alpha$.

Theorem~\ref{th:small_lambda} can be proved by using the \emph{dummy} balls idea used in the proof of Theorems~3.21 and~3.22 in~\cite{Janson04}: 
\begin{itemize}
\item Consider the continuous time urn with $d+1$ colours, such that the $d+1$ first colours evolve as the original $d$--colour process except that each time a ball splits, one ball of colour $d+1$ is added to the process. When a ball of colour $d+1$ splits, it splits into itself, adding no new balls in the process.
\item Apply Theorem~\ref{th:Janson_small} to this $d+1$--colour process (this process satisfies $\mathtt{(S)}$).
\item Go back to the original $d$--colour process by appropriate projection.
\end{itemize}
We do not develop the proof since no new idea is needed from there.

\subsection{Projections on large Jordan spaces}\label{sub:large}
Except for this digression on small eigenvalues, we are interested in the present paper in the behaviour of $U(n)$ along large Jordan blocks.
\textbf{We will from now on fix $E$ a Jordan subspace of $R$ associated to a large eigenvalue $\lambda$.} We denote by $\nu+1$ the size of its associated Jordan block (being also the dimension of $E$) and we denote by $v$ one eigenvector of $E$ associated to the eigenvalue $\lambda$.

\subsubsection{State of the art}
The asympotic behaviour of $U(n)$ projected onto the subspace $E$ is described by the following theorem:
\begin{theorem}[cf.~\cite{Pouyanne08}]
\label{th:limitDT}
Under $\mathtt{(B)}$, $\mathtt{(T)}$ and $\mathtt{(S)}$,
if $\frac12<\sigma<1$, then,
\[\lim_{n\to\infty}\frac{\pi_E(U(n))}{n^{\nicefrac{\lambda}{S}}\ln^{\nu} n} = \frac{1}{\nu !} W^{DT} v,\]
a.s. and in all $L^p$ ($p\geq 1$), where $\pi_E(U(n))$ is the projection of the composition vector at time $n$ onto $E$ (according to the Jordan decomposition of $R$).
\end{theorem}

\begin{remark}
Note that different choices for $v$ are possible, and that the random variable $W^{DT}$ depends on this choice. The random variable $W^{DT}$ should actually be denoted by $W^{DT}_{E,v}$ since it depends on the Jordan subspace $E$ fixed and on the choice of $v$. For clarity's sake, we will stick to the ambiguous but simpler notation $W^{DT}$; there is no ambiguity since $E$ and $v$ are fixed all along the present paper. Note that one could also choose to include the $1/\nu !$ into the definition of $W^{DT}$; we choose to leave it outside as done in~\cite{Pouyanne08}. 
\end{remark}

In continuous time, the composition vector projected on a large stable subspace satisfies
\begin{theorem}[see~\cite{Janson04}]
\label{th:limitCT}
Under $\mathtt{(B)}$, $\mathtt{(T)}$ and $\mathtt{(S)}$,
if $\frac12 < \sigma < 1$, then, almost surely and in all $L^p$ ($p\geq 1$),
\[\lim_{t\to+\infty}\frac{\pi_E(U(t))}{t^{\nu}\mathtt{e}^{\lambda t}} = \frac{1}{\nu !} W^{CT} v,\]
where $\pi_E$ is the projection on the large stable subspace $E$. 
Moreover, the random variable $W^{CT}$ admits moments of all orders.
\end{theorem}
\begin{remark}
Note that
Theorem~\ref{th:limitCT} is proven in~\cite{Janson04} under $\mathtt{(T_{-1})}$ and not $\mathtt{(T)}$; 
but Janson explains in his Remark~4.2 how to make it hold under $\mathtt{(T)}$.
\end{remark}

We are interested in the two random variables $W^{DT}$ and $W^{CT}$ defined in Theorems~\ref{th:limitDT} and~\ref{th:limitCT}. 
These random variables actually depend on the initial composition of the urn, denoted by $\bs{\alpha} = ^t\!(\alpha_1,\ldots,\alpha_d)$,
meaning that there are, for all $i\in\{1,\ldots,d\}$, $\alpha_i$ balls of type $i$ in the urn at time 0.
It is thus more rigorous to denote by $W^{DT}_{\bs{\alpha}}$ (resp. $W^{CT}_{\bs{\alpha}}$) the random variable associated to the initial composition $\bs \alpha$, emphasizing that we have to study two infinite families of random variables.

Connection~\eqref{eq:dcoul_connexion} implies connections between the random variables $W$ induced by the discrete and continuous processes. We need the following result to deduce them:
\begin{equation}\label{eq:dcoul_tau_n}
\lim_{n\to\infty}n\mathtt{e}^{-S\tau_n} = \xi,
\end{equation}
almost surely, where $\xi$ is a random variable with Gamma law of parameter $(\alpha_1+\ldots+ \alpha_d)/{S}$. 
This result is shown for two--colour urns in~\cite{CPS11}, and can be straightforwardly adapted to the present case, using the balanced hypothesis~$\mathtt{(B)}$. We do not develop this proof, which is very standard in the study of Yule processes (see for example~\cite[page~120]{AthreyaNey}).
We have
\[\frac{\pi_E\left(U^{CT}(\tau_n)\right)}{\tau_n^{\nu}\mathtt{e}^{\lambda\tau_n}} 
= \frac{\pi_E\left(U^{DT}(n)\right)}{n^{\nicefrac{\lambda}{S}}\ln^{\nu} n}\cdot \frac{n^{\nicefrac{\lambda}{S}}\ln^{\nu} n}{\tau_n^{\nu} \mathtt{e}^{\lambda \tau_n}}.\]
Moreover, Equation~\eqref{eq:dcoul_tau_n} implies $\frac{\ln n}{\tau_n}\to S$ when $n$ tends to $+\infty$,
\[\frac{n^{\nicefrac{\lambda}{S}}\ln^{\nu} n}{\tau_n^{\nu} \mathtt{e}^{\lambda \tau_n}} 
= \left(\frac{\ln n}{\tau_n}\right)^{\nu}(n\mathtt{e}^{-S\tau_n})^{\nicefrac{\lambda}{S}}\to S^{\nu}\xi^{\nicefrac{\lambda}{S}}.\]
Thus, for all initial composition $\bs \alpha$, we have (already mentioned in~\cite{Janson04}):
\begin{equation}
\label{eq:mcCT}
W^{CT}_{\bs \alpha} \EgalLoi S^{\nu} \xi^{\nicefrac{\lambda}{S}} W^{DT}_{\bs \alpha},
\end{equation}
where $\xi$ is a Gamma-distributed random variable with parameter $\frac{\alpha_1+\ldots+\alpha_d}{S}$,
and where $\xi$ and $W^{DT}$ are independent.

We also have that $(U^{DT}(n(t)))_{t\geq 0} = (U^{CT}(t))_{t\geq 0}$ almost surely, where $\sum_{i=1}^d \alpha_i + S n(t)$ is the total number of balls in the continuous time urn at time $t$. It implies that, for all initial composition $\bs \alpha$,
\begin{equation}
\label{eq:mcDT}
W^{DT}_{\bs \alpha} \EgalLoi S^{-\nu} \xi^{\nicefrac{-\lambda}{S}} W^{CT}_{\bs \alpha},
\end{equation}
where $\xi$ is a Gamma-distributed random variable with parameter $\frac{\alpha_1+\ldots+\alpha_d}{S}$ but where $\xi$ and $W^{CT}_{\bs \alpha}$ are \emph{not} independent, which can be verified via a covariance calculation.

\subsubsection{Statement of the main results}
The aim of the present paper is to gather information about $W^{DT}_{\bs \alpha}$ and $W^{CT}_{\bs \alpha}$.
We will prove the two following theorems, 
which happen to be the first results on the variables $W$ induced by a multi--colour urn.

\medskip
\begin{definition}\label{df:moment_det}
A complex-valued random variable~$Z$ is moment-determined if
for all random variable~$Y$,
\[\mathbb E\big[Y^p \bar Y^q\big] = \mathbb E\big[Z^p \bar Z^q\big]\quad (\forall p,q\geq 1)
\quad\Rightarrow\quad Y\EgalLoi Z.\]
\end{definition}

\begin{theorem}
\label{thm:dcoul_moments}
Under assumptions $\mathtt{(B)}$, $\mathtt{(T)}$ and $\mathtt{(S)}$,
for all initial composition~$\bs \alpha$,
\begin{enumerate}[(i)]
\item the random variable $W^{CT}_{\bs{\alpha}}$ is moment-determined.
\item the Laplace series of $W_{\bs{\alpha}}^{DT}$ converges on the whole plane, 
which implies that $W_{\bs{\alpha}}^{DT}$ is moment-determined.
\end{enumerate}
\end{theorem}
Note that it is an open problem to determine whether the Laplace transform of~$W^{CT}$ converges in a neighbourhood of~$0$.
A similar result is already proved in~\cite{CMP13} in the two--colour case 
under assumptions~$\mathtt{(B)}$, $\mathtt{(T_{-1})}$ and $\mathtt{(I)}$. 
Our proof is similar to the one developed there, although
more involved due to the higher dimension and to the less restrictive tenability assumption:
new higher--level arguments are needed here when simple calculations were sometimes enough in the two--colour case.

\medskip
\begin{theorem}
\label{th:density}
Under assumptions $\mathtt{(B)}$, $\mathtt{(T)}$ and $\mathtt{(I)}$ for all initial composition $\bs \alpha$:
\begin{itemize}
\item If $\lambda\in\mathbb{C}\setminus\mathbb{R}$, then the random variables $W_{\bs{\alpha}}^{CT}$ and $W_{\bs{\alpha}}^{DT}$ both admit a density on $\mathbb{C}$, and their support is the whole complex plane. 
\item If $\lambda\in\mathbb{R}$, then, the random variables $W_{\bs{\alpha}}^{CT}$ and $W_{\bs{\alpha}}^{DT}$ both admit a density on $\mathbb{R}$, and their support is the whole real line.
\end{itemize}
\end{theorem}
A similar result is proved in~\cite{CPS11} or in~\cite{CMP13} for two--colours urns 
with the tenability condition~$\mathtt{(T_{-1})}$
where the random variables $W^{CT}$ and $W^{DT}$ are real.
Note that in dimension~2, an urn with negative entries in its replacement matrix
admits no large Jordan block.
Since the random variables $W^{CT}$ and $W^{DT}$ can be non-real in the multi--colour case, 
the proof of Theorem~\ref{th:density} needs additional input;
for example, proving that the support of $W^{DT}$ is the whole complex plane 
becomes a non trivial step in the present article
whereas the fact that the support was the whole real line it was straightforward in the real case.

\section{Continuous time branching process -- smoothing system}\label{sec:CT_arbo}
In this section, we focus on the continuous--time process and show how to use the tree--like structure of the process. 
It is the first step of the proofs of our two main results Theorems~\ref{thm:dcoul_moments} and~\ref{th:density}.
As already mentioned, seeing urn schemes as branching processes is standard in the literature; 
the novelty of our approach is that we extend this analogy for urns that do not fall under the strong tenability condition~$\mathtt{(T_{-1})}$ but only under the weaker~$\mathtt{(T)}$.

First, we reduce the study to only $d$ initial compositions (instead of an infinite number), namely the initial compositions with a unique ball. Said differently, it is enough to study the random variables $W_{\bs{e_1}}, \ldots, W_{\bs{e_d}}$ where, for all $i\in\{1,\ldots, d\}$, $\bs{e_i}$ is the vector whose coordinates are all 0 except the $i^{\text{th}}$ which is 1 if $a_{i,i}\geq -1$ and $-a_{i,i}$ otherwise. We call $\bs{e_1}, \ldots, \bs{e_d}$ the {\bf atomic initial compositions} of the urn.
We then show, again using the tree-like structure of the process, that the random variables $W_{\bs{e_1}}, \ldots, W_{\bs{e_d}}$ 
satisfy a system of $d$ smoothing equations.

First introduce further notations:
For all $i\in\{1, \ldots, d\}$, let us denote
\begin{equation}\label{eq:tilde_alpha}
\tilde \alpha_i = \left\{\begin{array}{ll}
\alpha_i & \text{ if } a_{i,i}\geq -1,\\
& \\
\displaystyle \frac{\alpha_i}{-a_{i,i}} & \text{ otherwise},
\end{array}\right.
\end{equation}
and for all $i,c\in\{1, \ldots, d\}$, let us denote
\begin{equation}\label{eq:tilde_a}
\tilde a_{c,i} = \left\{\begin{array}{ll}
a_{c,i} & \text{ if } a_{i,i}\geq -1,\\
& \\
\displaystyle \frac{a_{c,i}}{-a_{i,i}} & \text{ otherwise}.
\end{array}\right.
\end{equation}
In view of Assumption~$\mathtt{(T)}$, for all $i,c\in\{1, \ldots, d\}$, $\tilde a_{c,i}$ and $\tilde \alpha_i$ are integers.

\begin{remark}
If we suppose further that $\mathtt{(T_{-1})}$ holds, then $\tilde \alpha_i = \alpha_i$ for all $i\in\{1, \ldots, d\}$ and $\tilde a_{i,j} = a_{i,j}$ for all $i,j\in\{1, \ldots, d\}$.
\end{remark}

\subsection{Decomposition}
To explain how to decompose the continuous time urn process, we will focus on an example, before generalising to any urn process that satisfies $\mathtt{(B)}$, $\mathtt{(S)}$ and $\mathtt{(T)}$. Assume for example that
\[R = \begin{pmatrix}6 & 2 & 0 \\ 5 & -2 & 5 \\ 0 & 2 & 6\end{pmatrix}.\]
One can verify that the urn process defined by $R$ satisfies $\mathtt{(B)}$, $\mathtt{(I)}$ and $\mathtt{(T)}$. 
Moreover, its eigenvalues are $8$, $6$ and $-4$. 
In particular, $6$ is a large eigenvalue which allows us to apply Section~\ref{sec:preliminaires} and define $W_{\bs \alpha}^{CT}$ through Theorem~\ref{th:limitCT} for any initial composition $\bs \alpha$. In the following, we will denote by $E_6$ the one-dimensional Jordan stable subspace associated to this eigenvalue $6$ and by $\pi_6$ the Jordan projection onto it.

\begin{figure}
\caption{Decomposition of the urn process in continuous time -- example.}
\label{fig:decomp}
\begin{center}
\includegraphics[width=.6\textwidth]{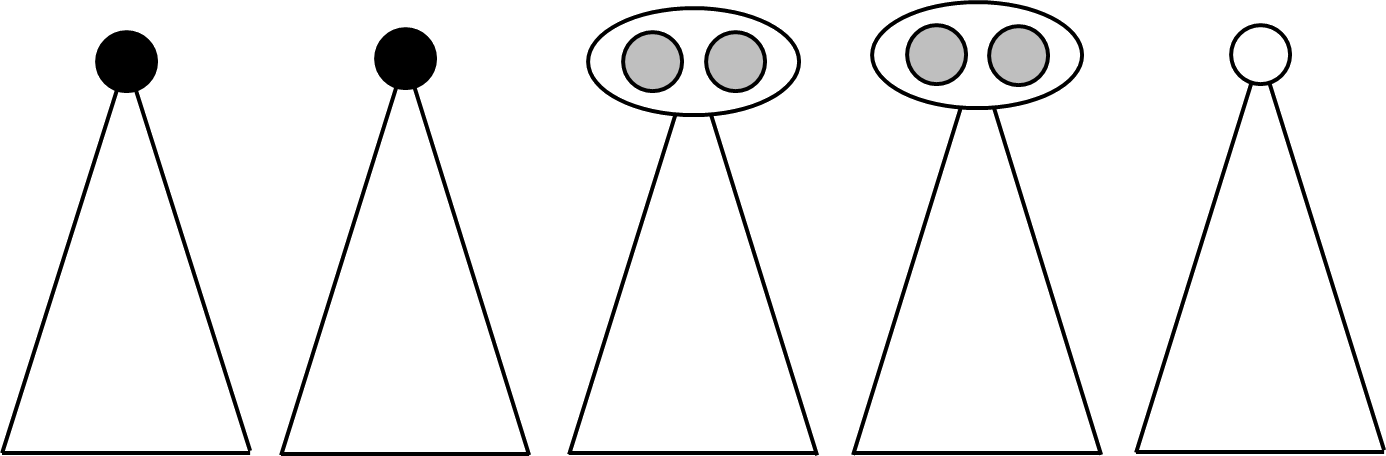}
\end{center}
\end{figure}

Note that we can decompose the multitype branching process as shown in Figure~\ref{fig:decomp}, which gives the following
\[U^{CT}_{(2, 4, 1)} \EgalLoi \sum_{k=1}^{2} U^{(k)}_{(1,0,0)} + \sum_{k=3}^4 U^{(k)}_{(0,2,0)} + \sum_{k=5}^5 U^{(k)}_{(0,0,1)},\]
where the $U^{(p)}$ are independent urn processes.

Let us quit the example and make the same reasoning as above in full generality under $\mathtt{(B)}$, $\mathtt{(S)}$ and $\mathtt{(T)}$.
Recall that, for all $c\in\{1, \ldots, d\}$, $\boldsymbol{e_c}$ has all its coordinates equal to zero except for the $c^{\text{th}}$, which is equal to $1$ if $a_{c,c}\geq 0$ and to $-a_{c,c}$ otherwise. We get
\[U^{CT}_{\boldsymbol{\alpha}}(t) \EgalLoi \sum_{c=1}^{d} \sum_{p=\beta_{c-1}+1}^{\beta_c} U_{\boldsymbol{e_c}}^{(p)}(t),\]
where $\beta_0 = 0$ and $\beta_c=\sum_{j\leq c} \tilde \alpha_j$, and where the $U_{\boldsymbol{e_c}}^{(p)}(t)$ are independent copies of $U^{CT}_{\boldsymbol{e_c}}(t)$, independent of each other.

Dividing this equality in law by $t^{\nu}\mathtt{e}^{\lambda t}$, projecting onto the fixed large Jordan subspace~$E$ via~$\pi_E$, 
and applying Theorem~\ref{th:limitCT} gives
\begin{proposition}[{already mentioned in~\cite[Remark~4.2]{Janson04}}]
\label{th:dcoul_decompCT}
For all replacement matrices $R$ and initial composition $\bs \alpha$ satisfying $\mathtt{(B)}$, $\mathtt{(S)}$ and $\mathtt{(T)}$,
\[W^{CT}_{\boldsymbol{\alpha}} \EgalLoi \sum_{c=1}^{d} \sum_{p=\beta_{c-1}+1}^{\beta_c} W_{\boldsymbol{e_c}}^{(p)},\]
where the $W_{\boldsymbol{e_c}}^{(p)}$ are independent copies of $W^{CT}_{\boldsymbol{e_c}}$, independent of each other and independent of $U$.
\end{proposition}

Proposition~\ref{th:dcoul_decompCT} allows to reduce the study to only $d$ random variables, namely $(W_{\bs{e_1}}^{CT},\ldots,W_{\bs{e_d}}^{CT})$ instead of having to study an infinite family of such variables. Any information gathered about those $d$ random variables will a priori give us some information about any $W_{\bs{\alpha}}^{CT}$.

\subsection{Dislocation}
In view of Proposition~\ref{th:dcoul_decompCT}, it is enough to focus on the $d$ atomic initial compositions $\bs{e_1}, \ldots, \bs{e_d}$. Recall that for all $i\in\{1, \ldots, d\}$, $\bs{e_i}$ is the vector whose all components are zero, except the $i^{\text{th}}$ which is equal to $1$ if $a_{i,i}\geq 0$ and to $-a_{i,i}$ otherwise. We will from now on denote by $\theta_i\in\{1, 2, \ldots\}$ the non zero component of $\bs{e_i}$:
\begin{equation}\label{eq:theta_i}
\theta_i = \begin{cases}
1 & \text{ if } a_{i,i}\geq 0\\
-a_{i,i} &\text{ otherwise}.
\end{cases}
\end{equation}

\begin{figure}
\caption{Dislocation of a continuous time urn process -- the different atomic initial compositions and their composition after the first drawing.}
\label{fig:disloc}
\begin{center}
\includegraphics[width=.5\textwidth]{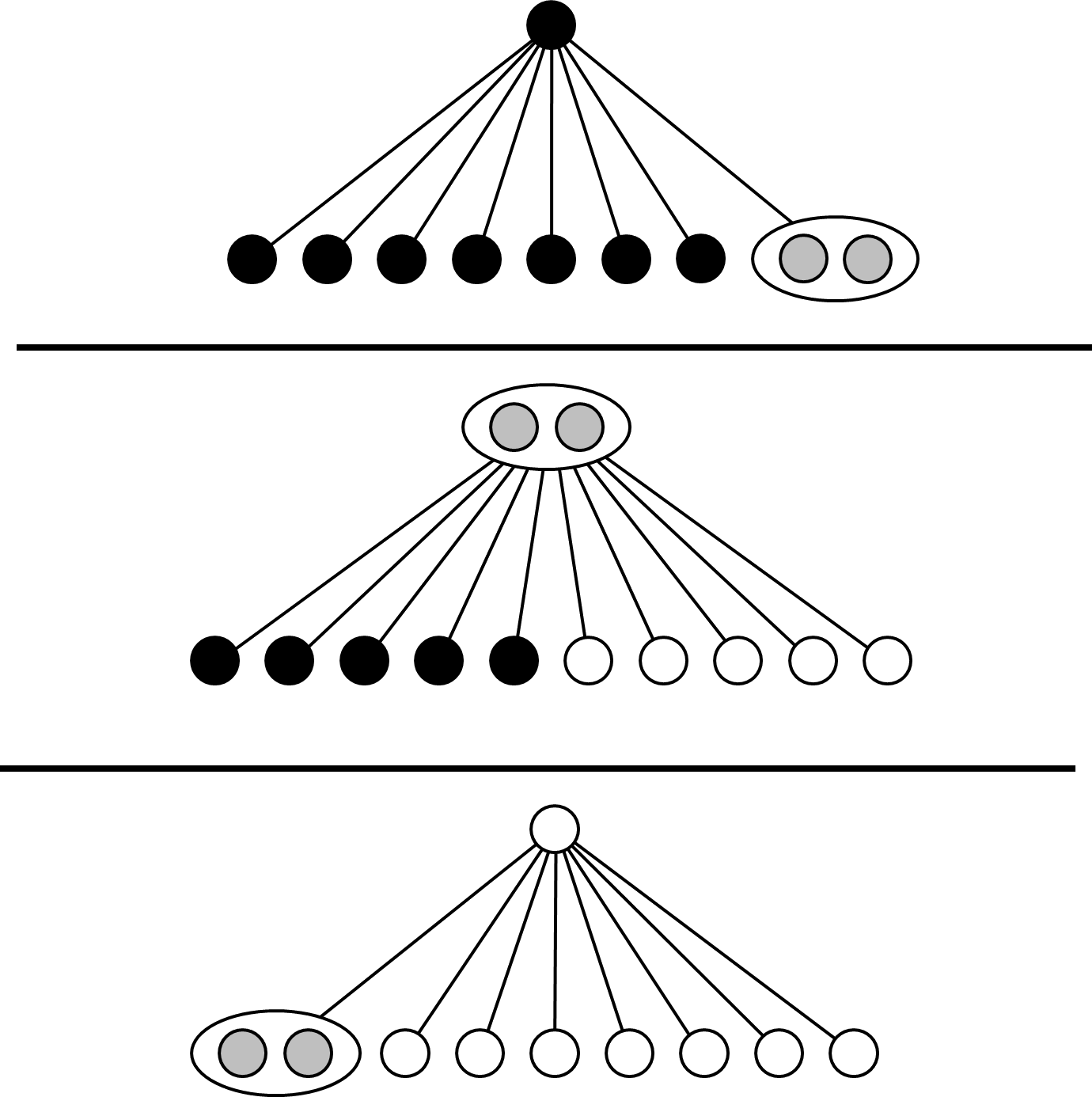}
\end{center}
\end{figure}

Let us again study first the particular example given by
\[R = \begin{pmatrix}6 & 2 & 0 \\ 5 & -2 & 5 \\ 0 & 2 & 6\end{pmatrix}.\]
The three atomic initial compositions are given by $\bs{e_1} = (1, 0, 0)$, $\bs{e_2} = (0, 2, 0)$ and $\bs{e_3} = (0, 0, 1)$, since $a_{1,1}, a_{3,3}\geq 0$ and $a_{2,2} = -2$. In all three cases, the first step is deterministic: we know the colour of the first ball to be drawn and we therefore know what is the composition of the urn after the first split time (cf.~Figure~\ref{fig:disloc}). We therefore have that,
\[U^{CT}_{\bs{e_1}} \EgalLoi \sum_{k=1}^{7} U^{(k)}_{\bs{e_1}}(t-\tau^{(1)}) + \sum_{k=8}^8 U^{(k)}_{\bs{e_2}}(t-\tau^{(1)}),\]
\[U^{CT}_{\bs{e_2}} \EgalLoi \sum_{k=1}^{5} U^{(k)}_{\bs{e_1}}(t-\tau^{(2)}) + \sum_{k=6}^{10} U^{(k)}_{\bs{e_3}}(t-\tau^{(2)}),\]
\[U^{CT}_{\bs{e_3}} \EgalLoi \sum_{k=1}^{1} U^{(k)}_{\bs{e_2}}(t-\tau^{(3)}) + \sum_{k=2}^8 U^{(k)}_{\bs{e_3}}(t-\tau^{(3)}),\]
where the $U^{(k)}$ are independent continuous time urn processes with replacement matrix $R$, where $\tau^{(1)}$, $\tau^{(2)}$ and $\tau^{(3)}$ are independent random variables exponentially distributed of respective parameters $1$, $2$ and $1$. The random variables $\tau^{(1)}$, $\tau^{(2)}$ and $\tau^{(3)}$ are the first split times that occur in urns of respective initial compositions $\bs{e_1}$, $\bs{e_2}$ and~$\bs{e_3}$.

The same reasoning in full generality, for all replacement matrices $R$ satisfying $\mathtt{(B)}$, $\mathtt{(S)}$ and $\mathtt{(T)}$ gives that, for all $c\in\{1, \ldots, d\}$,
\[U^{CT}_{\bs{e_c}} \EgalLoi \sum_{i=1}^d \sum_{k=\gamma^{(c)}_{i-1}+1}^{\gamma^{(c)}_i} U^{(k)}_{\bs{e_i}}(t-\tau^{(c)}),\]
where $\gamma^{(c)}_0 = 0$ and $\gamma_i^{(c)} = \sum_{j\leq i} \tilde a_{c,i} + \delta_{c,i}$, 
$\tau^{(c)}$ is an exponentially distributed random variable of parameter $\theta_c$ (which is a positive integer), 
and the $U^{(k)}$ are independent continuous time urn processes with replacement matrix $R$.
We recall that $E$ is a fixed Jordan space of dimension $\nu+1$ associated to a large eigenvalue $\lambda$.
Dividing the previous equality in law by $t^{\nu}\mathtt{e}^{\lambda t}$, 
projecting onto $E$ via $\pi_E$, and applying Theorem~\ref{th:limitCT} gives

\begin{proposition}[{already mentioned in~\cite[Theorem~3.9]{Janson04}}]\label{thm:dcoul_dislocCT}
Under assumptions $\mathtt{(B)}$, $\mathtt{(S)}$ and $\mathtt{(T)}$,
for all $c\in\{1,\ldots,d\}$,
\begin{equation}
\label{eq:dcoul_sysCT}
W^{CT}_{\boldsymbol{e_c}} \EgalLoi U^{\nicefrac{\lambda}{\theta_c}} \sum_{i=1}^d \sum_{k=\gamma_{i-1}^{(c)}+1}^{\gamma_i^{(c)}} W^{(k)}_{\boldsymbol{e_i}},
\end{equation}
where $\gamma^{(c)}_0 = 0$, $\gamma_i^{(c)} = \sum_{j\leq i}(\tilde a_{c,j}+\delta_{c,j})$ (using Kronecker's notation $\delta_{c,j}=1$ if $c=j$, 0 otherwise), where $U$ is a uniform random variable on $[0,1]$, and where the $W^{(k)}_{\boldsymbol{e_i}}$ are independent copies of $W^{CT}_{\boldsymbol{e_i}}$, independent of each other and of $U$.
\end{proposition}

\begin{remark}
If we assume $\mathtt{(T_{-1})}$ instead of $\mathtt{(T)}$ in the result above, we get
\begin{equation*}
W^{CT}_{\boldsymbol{e_c}} 
\EgalLoi U^{\lambda} \sum_{i=1}^d \sum_{k=\gamma_{i-1}^{(c)}+1}^{\gamma_i^{(c)}} W^{(k)}_{\boldsymbol{e_i}},
\end{equation*}
where $\gamma^{(c)}_0 = 0$, $\gamma_i^{(c)} = \sum_{j\leq i}(a_{c,j}+\delta_{c,j})$, where $U$ is a uniform random variable on $[0,1]$, and where the $W^{(k)}_{\boldsymbol{e_i}}$ are independent copies of $W^{CT}_{\boldsymbol{e_i}}$, independent of each other and of $U$.
\end{remark}

\begin{remark}
One can prove that the solution of System~\eqref{eq:dcoul_sysCT} 
is unique at fixed mean and under a condition of finite variance. 
A very similar proof in done in~\cite[{Proof of Theorem~3.9(iii), page~232--233}]{Janson04}.
\end{remark}

\section{Moments -- Proof of Theorem~\ref{thm:dcoul_moments}}\label{sec:dcoul_moments}
This section is devoted to the proof of Theorem~\ref{thm:dcoul_moments}.
We study here the moments of the random variables $(W^{DT}_{\bs{e_i}})_{i\in\{1,\ldots,d\}}$ 
and $(W^{CT}_{\bs{e_i}})_{i\in\{1,\ldots,d\}}$. 
First note that the convergence in all $L^p$ ($p\geq 1$) stated in Theorems~\ref{th:limitDT} and~\ref{th:limitCT} 
ensures us that those random variables admit moments of all orders.
The first step in the proof is the following lemma, which concerns the continuous time process: 
Equation~\eqref{eq:mcDT} will then allow us to infer results about the discrete time process.
\begin{lemma}
\label{lem:dcoul_borne}
Let $(X_1,\ldots, X_d)$ be a solution of System~\eqref{eq:dcoul_sysCT} with moments of all orders. Then the sequences $\displaystyle\left(\frac{\mathbb{E} |X_i|^p}{p!\ln^p p}\right)^{\frac1p}$, for all $i\in\{1,\ldots,d\}$, are bounded.
\end{lemma}

\begin{proof}
Let $(X_1,\ldots, X_d)$ be a solution of System~\eqref{eq:dcoul_sysCT}, let $\varphi(p):= \ln^p (p+2)$ and let, for all $i\in\{1,\ldots, d\}$,
\[u_p^{(i)}:=\frac{\mathbb{E} |X_i|^p}{p!\varphi(p)}.\]
Let us prove by induction on $p\geq 1$ that, for all $i\in\{1,\ldots, d\}$, the sequence $\left(\frac{\mathbb{E} |X_i|^p}{p!\varphi (p)}\right)^{\frac 1p}$ is bounded.
Raise the equations of System~\eqref{eq:dcoul_sysCT} to the power~$p$ and apply the multinomial formula. 
Since $\mathbb{E}|U^{\mu p}| = \frac1{p\RRe\mu+1}$ for all $\mu\in\mathbb{C}$ (with $U$ uniformly distributed on $[0,1]$), for all $c\in\{1,\ldots, d\}$,
\begin{align*}
\mathbb{E}|X_c|^p 
\leq & 
\frac 1{p\frac{\RRe \lambda}{\theta_c} +1} \bigg(\sum_{i=1}^d (\tilde a_{c,i}+\delta_{c,i}) \mathbb{E} |X_i|^p \\
&+ \summ{p_1+ \dots + p_{\gamma_d^{(c)}} = p}{ p_j\leq p-1} \frac{p!}{\prod_{1\le j\le\gamma_d^{(c)}}p_j!}
\prod_{i=1}^d \prod_{k=\gamma_{i-1}^{(c)}+1}^{\gamma_i^{(c)}}\mathbb{E}|X_i|^{p_{k}}\bigg),
\end{align*}
which means
\[
p\frac{\RRe \lambda}{\theta_c} \mathbb{E}|X_c|^p 
\leq  \sum_{i=1}^d \tilde a_{c,i} \mathbb{E}|X_i|^p 
+ \summ{p_1+ \ldots + p_{\gamma_d^{(c)}} = p}{ p_j\leq p-1} \frac{p!}{\prod_{1\le j\le\gamma_d^{(c)}}p_j!}
\prod_{i=1}^d \prod_{k=\gamma_{i-1}^{(c)}+1}^{\gamma_i^{(c)}}\mathbb{E}|X_i|^{p_{k}}.
\]
It implies that, for all $c\in\{1,\ldots, d\}$,
\begin{equation}
\label{eq:dcoul_inegaliteMoments}
\frac{p\RRe \lambda}{\theta_c} u^{(c)}_p \leq
\sum_{i=1}^d \tilde a_{c,i} u_p^{(i)} + 
\summ{p_1+ \ldots + p_{\gamma_d^{(c)}} = p}{ p_j\leq p-1} \frac{\prod_{1\le j\le \gamma_d^{(c)}}\varphi(p_j)}{\varphi(p)} \prod_{i=1}^d \prod_{k=\gamma_{i-1}^{(c)}+1}^{\gamma_i^{(c)}} u^{(i)}_{p_k}.
\end{equation}
Let
\begin{equation}
\label{dcoul_Phi}
\Phi_c(p):= \summ{p_1+ \dots + p_{\gamma_d^{(c)}} = p}{ p_j\leq p-1} \   \frac{\prod_{1\le j\le \gamma_d^{(c)}}\varphi(p_j)}{\varphi(p)}.
\end{equation}
A slight generalisation of~\cite[Lemma~1]{CMP13} ensures that 
\begin{equation}\label{eq:phi_c}
\Phi_c(p)\leq \left(1+8\ln(p+2)\right)^{\gamma_d^{(c)}},
\end{equation}
for all $p\geq 2$, as soon as $\gamma_d^{(c)}\geq 1$; 
recall that $\gamma_d^{(c)} = \sum_{i=1}^d (\tilde a_{c,i}+\delta_{c,i}) \geq 1$.

Denote by $\Delta_p$ the determinant of $p\RRe\lambda \Theta-\tilde R$ where $\Theta = \mathtt{diag}\{\theta_1^{-1}, \ldots, \theta_d^{-1}\}$ and $\tilde R = (\tilde a_{i,j})_{1\leq i,j\leq d}$. 
This determinant is non-zero for all $p\geq 2$ since 
\[\Delta_p = \mathtt{det}(p\RRe\lambda \Theta-\tilde R) = \bigg(\prod_{i=1}^d \theta_i^{-1}\bigg) \mathtt{det}(p\RRe\lambda I_d - R)\] 
and since $p\RRe\lambda > S$ (as, by assumption, $\RRe\lambda>S/2$). 
Note that, under $\mathtt{(T_{-1})}$, we have $\Theta= I_d$.

Recall that $\theta_i \geq 1$ for all $i\in\{1, \ldots, d\}$ (see Equation~\eqref{eq:theta_i}), and note that
\[\frac{\|p\RRe\lambda \Theta - \tilde R\|_{\infty}}{p \RRe\lambda}\to \frac{1}{\min_{1\le i\le d}\theta_i}\] 
when $p$ tends to $+\infty$, implying that there exists $p_0\geq 1$ such that, for all $p\geq p_0$,
\begin{equation}\label{eq:rho}
\frac{\|p\RRe\lambda \Theta - \tilde R\|_{\infty}}{p\RRe\lambda}\leq \frac{2}{\min \theta_i}=: \rho.
\end{equation}
In addition, let us denote by $\Delta_p(j,i)$ the determinant of $p\RRe\lambda \Theta - \tilde R$ 
in which the $i^{\text{th}}$ column and the $j^{\text{th}}$ line have been removed. 
For all $1\leq i,j\leq d$, the polynomial $\Delta_p(j,i)$ has degree at most $d-1$ in $p$, 
which implies
\[\sup_{1\leq i,j\leq d} \frac{|\Delta_p(j,i)|}{|\Delta_p|} = \mathcal{O}\left(\frac1p\right),\]
when $p$ goes to infinity,
and there exists a constant $\eta>0$ and an integer $p_1\geq p_0$ such that, for all $p\geq p_1$,
\[\sup_{1\leq i,j\leq d} \frac{|\Delta_p(j,i)|}{|\Delta_p|} \leq \frac{\eta}{p}.\]
Finally, let us denote by $\Delta_p(c)$ the determinant of the matrix $p\RRe\lambda \Theta - \tilde R$ 
in which the $c^{\text{th}}$ column has been replaced by a column of 1. 
We know that $\Delta_p$ has degree $d$ in $p$ whereas $\Delta_p(c)$ is a polynomial with degree at most $d-1$ in $p$. 
It implies that there exists an integer $p_2\geq p_1$ such that, for all $p\geq p_2$, for all $c\in\{1,\ldots, d\}$,
\begin{equation}\label{eq:num}
\frac{\Delta_p(c)}{\Delta_p} \left(1+8 \ln(p+2)\right)^{\gamma_d^{(c)}} \leq \frac{1}{\rho\eta d^2\RRe\lambda},
\end{equation}
where the choice of the right-hand side constant will become clear later on.

Let us define
\[A:= \max\{(u_q^{(i)})^{\frac 1q}, 1\leq q\leq p_2, 1\leq i\leq d\},\]
and prove by induction on $p\geq p_2$ that, for all $q\leq p$ and $c\in\{1,\ldots,d\}$, $(u_q^{(c)})^{\frac 1q} \leq A$.
Fix $p>p_2$ and assume that the induction hypothesis is true for $p-1$; then
Equations~\eqref{eq:phi_c} and~\eqref{eq:dcoul_inegaliteMoments} imply
\[
p\RRe\lambda \Theta u_p^{(c)} \leq \sum_{i=1}^d \tilde a_{c,i} u_p^{(i)} + A^p\Phi_c(p).
\]
Let $(v_1, \ldots, v_d)$ be the solution of the system
\[p\RRe\lambda \Theta v_c = \sum_{i=1}^d \tilde a_{c,i} v_i + A^p\Phi_c(p).\]
We thus have that $v_c \geq u_p^{(c)}$ for all $1\leq c\leq d$.
Solving this Cramér system, we get, in view of Equations~\eqref{eq:phi_c} and~\eqref{eq:num},
\[v_c = A^p \Phi_c(p) \frac{\Delta_p(c)}{\Delta_p}\leq \frac{A^p}{\rho\eta d^2\RRe\lambda}.\]
For all $p\geq p_2$, using the fact that for all $d$-dimensional matrix $M$ and all vector $\bs v$, 
$\|M\bs v\|_{\infty}\leq d\|M\|_{\infty}\|\bs v\|_{\infty}$, we get
\[(p\RRe\lambda \Theta - \tilde R)\bs{u}^{(p)} 
\leq (p\RRe\lambda \Theta - \tilde R)\bs{v} 
\leq \|p\RRe\lambda \Theta -\tilde R\|_{\infty}  \frac{A^p}{\rho\eta d\RRe\lambda}  \bs\omega,\]
where $\bs u^{(p)}$ and $\bs v$ denote the vectors of respective coordinates $(u_p^{(i)})_{1\leq i\leq d}$ and $(v_i)_{1\leq i\leq d}$, where $\bs \omega$ is the vector whose all coordinates are equal to 1, and where the sign $\leq$ between two vectors is to be read coordinate by coordinate.
In particular, we have
\[\|(p\RRe\lambda \Theta -\tilde R)\bs{u}^{(p)}\|_{\infty} 
\leq A^p\frac{\|p\RRe\lambda \Theta -\tilde R\|_{\infty}}{\rho\eta d\RRe\lambda}
\leq A^p\frac{p}{\eta d},\]
where we have used Equation~\eqref{eq:rho}.
Let us denote $M = p\RRe\lambda \Theta - \tilde R$. The coefficients of $M^{-1}$ are given by
\[(M^{-1})_{i,j} = (-1)^{i+j}\frac{\Delta_p(j,i)}{\Delta_p},\]
where $\Delta_p$ is the determinant of $M$, and $\Delta_p(j,i)$ is the determinant of the matrix $M$ in which the $j^{\text{th}}$ line and the $i^{\text{th}}$ column have been removed.
By definition of $p_2$, for all $p\geq p_2$,
\[\|M^{-1}\|_{\infty} = \sup_{1\leq i,j\leq d} |(M^{-1})_{i,j}| \leq \frac{\eta}{p},\]
which implies that, for all $p\geq p_2$,
\[\|\bs{u}^{(p)}\|_{\infty} 
\leq \| M^{-1}\|_{\infty}  A^p \frac{p}{\eta} 
\leq A^p.\]
Finally, for all $c\in\{1,\ldots, d\}$,
$u_c^{(p)} \leq A^p$,
which concludes the proof.
\end{proof}

\begin{proof}[Proof of Theorem~\ref{thm:dcoul_moments}]
\hspace{-4pt}{\it (i)}\hspace{-2pt} To prove that a real-valued random variable~$X$ is moment-determined, 
one can show that it satisfies Carleman's criterion:
\[\sum_{k=1}^{\infty} \mathbb E\big[X^{2k}\big]^{-\nicefrac1{2k}} = \infty.\]
For a complex-valued random variable $Z$, one can for example apply~\cite[Theorem~10.3]{JK}, which states that,
if all moments of~$|Z|$ are finite and if $|Z|$ satisfies Carleman's criterion, 
then~$Z$ is moment-determined in the sense of Definition~\ref{df:moment_det}.
Lemma~\ref{lem:dcoul_borne} implies that, for all $i\in\{1,\ldots, d\}$, 
the random variable $|W^{CT}_{\bs {e_i}}|$, 
which admits moments of all orders in view of Theorem~\ref{th:limitCT} 
satisfies Carleman's criterion. Therefore, $W^{CT}_{\bs {e_i}}$ is moment-determined. 
Proposition~\ref{th:dcoul_decompCT} eventually allows us to generalise 
Lemma~\ref{lem:dcoul_borne} to any initial composition: 
For all initial composition $\bs\alpha$, $|W^{CT}_{\bs\alpha}|$ also satisfies the Carleman's criterion.

{\it (ii)} By Lemma~\ref{lem:dcoul_borne}, we have the following inequality: there exists a constant $C$ such that, for all integer $p\geq 2$, for any initial composition,
\[\frac{\mathbb{E}|W^{CT}|^p}{p!}\leq C^p \ln^p p.\]
It implies, via Equation~\eqref{eq:mcCT}, there exists a constant $D$ such that, for all integer $p\geq 2$,
\[\frac{\mathbb{E}|W^{DT}|^p}{p!}\leq D^p \frac{\ln^p p}{\Gamma\Big(\frac{p\RRe\lambda +1}{S}\Big)},\]
where we recall that $\RRe\lambda>S/2>0$.
This implies that the Laplace series of $|W^{DT}|$ has an infinite radius of convergence.
Since $\langle t, W^{DT}\rangle \leq 2|t||W^{DT}|$,
this implies that the Laplace transform of $W^{DT}$ converges on the whole complex plane.
\end{proof}

Theorem~\ref{thm:dcoul_moments} gives an upper bound for the moments of $W^{CT}$ and $W^{DT}$; note that no lower bound is known up to now.

\section{Discrete time urn process -- Smoothing system}\label{sec:discrete}

\subsection{Smoothing system in discrete time}
This subsection is devoted to deduce from Sections~\ref{sec:CT_arbo} and~\ref{sec:dcoul_moments} that the random variable $(W^{DT}_{\bs{e_1}}, \ldots, W^{DT}_{\bs{e_d}})$ is a solution of a smoothing system:
\begin{proposition}\label{thm:dcoul_dislocDT}
Under assumptions $\mathtt{(B)}$, $\mathtt{(T)}$ and $\mathtt{(S)}$,
for every colour $1\leq c\leq d$,
\begin{equation}
\label{eq:dcoul_sysDT}
W^{DT}_{\boldsymbol{e_c}} 
\EgalLoi 
\sum_{i=1}^d \sum_{k=\gamma_{i-1}^{(c)}+1}^{\gamma_i^{(c)}} \left(V^{(c)}_k\right)^{\nicefrac{\lambda}{S}} W_{\boldsymbol{e_i}}^{(k)},
\end{equation}
where $\gamma_0 = 0$ and $\gamma_i^{(c)} = \sum_{j=1}^i (\tilde a_{c,j}+\delta_{c,j})$ for all $i\in\{1,\ldots,d\}$; 
where the $W_{\boldsymbol{e_i}}^{(k)}$ are independent copies of $W^{DT}_{\boldsymbol{e_i}}$, independent of each other; 
and where $V^{(c)}=(V_1,\ldots,V_{\gamma^{(c)}_d})$ is a Dirichlet-distributed random vector independent of the $W$ and of parameter $\bs \pi$, given by
\[\pi_k = \nicefrac{\theta_i}{S} \quad\text{ if }\quad \gamma_{i-1}^{(c)} < k \leq \gamma_i^{(c)}.\]
\end{proposition}

\begin{remark}
If we assume $\mathtt{(T_{-1})}$ instead of $\mathtt{(T)}$ in the theorem above, we obtain the following system:
\begin{equation*}
W^{DT}_{\boldsymbol{e_c}} 
\EgalLoi 
\sum_{i=1}^d \sum_{k=\gamma_{i-1}^{(c)}+1}^{\gamma_i^{(c)}} V_k^{\nicefrac{\lambda}{S}} W_{\boldsymbol{e_i}}^{(k)},
\end{equation*}
where $\gamma_0 = 0$ and $\gamma_i^{(c)} = \sum_{j=1}^i (a_{c,j}+\delta_{c,j})$ for all $i\in\{1,\ldots,d\}$; where the $W_{\boldsymbol{e_i}}^{(k)}$ are independent copies of $W^{DT}_{\boldsymbol{e_i}}$, independent of each other; and where $V=(V_1,\ldots,V_{S+1})$ is a Dirichlet-distributed random vector of parameter $\left(\frac1{S},\ldots,\frac1{S}\right)$, independent of the $W$.
\end{remark}

\begin{proof}
We will present two proof for this statement: the first one is a moment proof, developed hereafter, the second one is the classical one based on the branching property of the urn process, we will detail this other proof in Subsection~\ref{sub:alt_proof}.

Let us prove that, for all $p, q\geq 1$,
\[\mathbb{E} \Big[\big(W^{DT}_{\boldsymbol{e_c}}\big)^p \big(\overline W^{DT}_{\boldsymbol{e_c}}\big)^q\Big]
= \mathbb{E} 
\bigg[\bigg(\sum_{i=1}^d \sum_{k=\gamma_{i-1}^{(c)}+1}^{\gamma_i^{(c)}} V_k^{\nicefrac{\lambda}{S}} W_{\boldsymbol{e_i}}^{(k)}
\bigg)^{\!p}
\bigg(\sum_{i=1}^d \sum_{k=\gamma_{i-1}^{(c)}+1}^{\gamma_i^{(c)}} V_k^{\nicefrac{\bar\lambda}{S}} \overline W_{\boldsymbol{e_i}}^{(k)}
\bigg)^{\!q}\bigg].\]
Since, $W^{DT}_{\bs{e_c}}$ is moment-determined in view of Theorem~\ref{thm:dcoul_moments}, this will conclude the proof.
Let us use Connection~\eqref{eq:dcoul_connexion}, which gives, for all $c\in\{1, \ldots, d\}$,
\begin{equation}\label{eq:moments_connexion}
\mathbb{E} \big[\big(W^{CT}_{\boldsymbol{e_c}}\big)^p \big(\overline W^{CT}_{\boldsymbol{e_c}}\big)^q\big] 
=  S^{(p+q)\nu} \mathbb{E}\Big[\xi_c^{\frac{p\lambda+q\bar\lambda}{S}}\Big] 
\, \mathbb{E} \big[\big(W^{DT}_{\boldsymbol{e_c}}\big)^p \big(\overline W^{DT}_{\boldsymbol{e_c}}\big)^q\big],
\end{equation}
where $\xi_c$ is a Gamma-distributed random variable, of parameter $\nicefrac{\theta_c}{S}$.

Let $V_c$ be a Beta-distributed random variable of parameter $\left(\nicefrac{\theta_c}{S},1\right)$. 
Note that if $U$ is uniformly distributed on $[0,1]$, we have
\[V_c \EgalLoi U^{\nicefrac{S}{\theta_c}}.\]
Let $(\zeta_{i,k})_{1\leq i\leq d; k\geq 1}$ be a sequence of independent, 
Gamma-distributed random variables of parameter $\nicefrac{\theta_c}{S}$.
The random variable
\[\zeta_c = V_c \sum_{i=1}^d\sum_{k=\gamma_{i-1}^{(c)}}^{\gamma_i^{(c)}} \zeta_{i,k}\]
is Gamma-distributed with parameter $\nicefrac{\theta_c}{S}$ (it can be verified by calculating its moments).
Finally, let
\[V^{(c)}_{i,k} = \frac{\zeta_{i,k}}{\sum_{i=1}^d\sum_{k=\gamma_{i-1}^{(c)}}^{\gamma_i^{(c)}} \zeta_{i,k}}.\]
Then (see for example~\cite[Lemma~2.2]{Bertoin06}), the $\gamma_d^{(c)}$-dimensional random vector $V^{(c)}$ whose coordinates are given by
\[V^{(c)}_k = V^{(c)}_{i,k} \quad \text{ if } \gamma_{i-1}^{(c)} < k \leq \gamma_i^{(c)}\]
is Dirichlet-distributed of parameter $\bs \pi = (\pi_1, \ldots, \pi_{\gamma_d^{(c)}})$, where
\[\pi_k = \nicefrac{\theta_i}{S} \quad\text{ if }\quad \gamma_{i-1}^{(c)} < k \leq \gamma_i^{(c)},\]
and independent from $\zeta_c$. 
Thus, System~\eqref{eq:dcoul_sysCT} together with Equation~\eqref{eq:moments_connexion}, gives that, for all $c\in\{1, \ldots, d\}$,
\begin{align}
&S^{(p+q)\nu} \mathbb{E}\Big[\xi_c^{\frac{p\lambda+q\bar\lambda}{S}}\Big] 
\, \mathbb{E} \big[\big(W^{DT}_{\boldsymbol{e_c}}\big)^{\!p} \big(\overline W^{DT}_{\boldsymbol{e_c}}\big)^{\!q}\big]\\
&= \mathbb{E}\bigg[\bigg(V_c^{\nicefrac{\lambda}{S}} \sum_{i=1}^d \sum_{k=\gamma^{(c)}_{i-1}+1}^{\gamma_i^{(c)}} 
W_{\bs{e_i}}^{CT,(k)}\bigg)^{\!\!p}
\bigg(V_c^{\nicefrac{\bar \lambda}{S}} \sum_{i=1}^d \sum_{k=\gamma^{(c)}_{i-1}+1}^{\gamma_i^{(c)}} 
\overline W_{\bs{e_i}}^{CT,(k)}\bigg)^{\!\!q}\bigg]\nonumber\\
&= 
\mathbb{E}\bigg[
\bigg(\sum_{i=1}^d\sum_{k=\gamma^{(c)}_{i-1}+1}^{\gamma_i^{(c)}} 
V_c^{\nicefrac{\lambda}{S}} S^{\nu} \zeta_{i,k}^{\nicefrac{\lambda}{S}} W_{\bs{e_i}}^{DT,(k)}\bigg)^{\!\!p}
\bigg(\sum_{i=1}^d\sum_{k=\gamma^{(c)}_{i-1}+1}^{\gamma_i^{(c)}} 
V_c^{\nicefrac{\bar\lambda}{S}} S^{\nu} \zeta_{i,k}^{\nicefrac{\bar\lambda}{S}} \overline W_{\bs{e_i}}^{DT,(k)}\bigg)^{\!\!q}\bigg]\nonumber\\
&= S^{(p+q)\nu}\;
\mathbb{E}\bigg[
\bigg(\sum_{i=1}^d\sum_{k=\gamma^{(c)}_{i-1}+1}^{\gamma_i^{(c)}} \!\!\!
\big(\zeta_{i,k}V_c\big)^{\!\nicefrac{\lambda}{S}} W_{\bs{e_i}}^{DT,(k)}\bigg)^{\!\!p}
\bigg(\sum_{i=1}^d\sum_{k=\gamma^{(c)}_{i-1}+1}^{\gamma_i^{(c)}} \!\!\!
\big(\zeta_{i,k}V_c\big)^{\!\nicefrac{\bar\lambda}{S}} \overline W_{\bs{e_i}}^{DT,(k)}\bigg)^{\!\!q}\bigg].\label{eq:first_step}
\end{align}
We also have that, by independence of $V^{(c)}$ and $\zeta_c$,
\begin{align}
&\mathbb{E}\Big[\zeta_c^{\frac{p\lambda+q\bar\lambda}{S}}\Big]
\mathbb{E}\bigg[\bigg(\sum_{i=1}^d\sum_{k=\gamma^{(c)}_{i-1}+1}^{\gamma_i^{(c)}} \!\!\!
\left(V^{(c)}_{i,k}\right)^{\!\nicefrac{\lambda}{S}} W_{\bs{e_i}}^{DT,(k)}\bigg)^{\!\!p}
\bigg(\sum_{i=1}^d\sum_{k=\gamma^{(c)}_{i-1}+1}^{\gamma_i^{(c)}} \!\!\!
\left(V^{(c)}_{i,k}\right)^{\!\nicefrac{\bar \lambda}{S}} \overline W_{\bs{e_i}}^{DT,(k)}\bigg)^{\!\!q}\bigg]\nonumber\\
&= 
\mathbb{E}\bigg[\bigg(\sum_{i=1}^d\sum_{k=\gamma^{(c)}_{i-1}+1}^{\gamma_i^{(c)}} \!\!\!
\left(\zeta_c V^{(c)}_{i,k}\right)^{\!\nicefrac{\lambda}{S}} W_{\bs{e_i}}^{DT,(k)}\bigg)^{\!\!p}
\bigg(\sum_{i=1}^d\sum_{k=\gamma^{(c)}_{i-1}+1}^{\gamma_i^{(c)}} \!\!\!
\left(\zeta_c V^{(c)}_{i,k}\right)^{\!\nicefrac{\bar \lambda}{S}} \overline W_{\bs{e_i}}^{DT,(k)}\bigg)^{\!\!q}\bigg]\nonumber\\
&
\mathbb{E}\bigg[\bigg(\sum_{i=1}^d\sum_{k=\gamma^{(c)}_{i-1}+1}^{\gamma_i^{(c)}} \!\!\!
\left(\zeta_{i,k}V_c\right)^{\!\nicefrac{\lambda}{S}} W_{\bs{e_i}}^{DT,(k)}\bigg)^{\!\!p}
\bigg(\sum_{i=1}^d\sum_{k=\gamma^{(c)}_{i-1}+1}^{\gamma_i^{(c)}} \!\!\!
\left(\zeta_{i,k}V_c\right)^{\!\nicefrac{\bar\lambda}{S}} \overline W_{\bs{e_i}}^{DT,(k)}\bigg)^{\!\!q}\bigg].\label{eq:second_step}
\end{align}
The result follows from Equations~\eqref{eq:first_step} and~\eqref{eq:second_step}, using the fact that~$W^{DT}$ is moment-determined.
\end{proof}

\subsection{Unicity}\label{sec:dcoul_contraction}
The main goal of this section is to prove that the solution of System~\eqref{eq:dcoul_sysDT} is unique. We therefore use the so-called contraction method. This method, presented for example in Neininger-R{\"u}schendorf's survey~\cite{NeiRusSurvey} consists in applying the Banach fixed point theorem in an appropriate complete Banach space. It has already been used in a P{\'o}lya urn context in the literature. In~\cite{KN13} the contraction method is used as a new approach to prove an equivalent of Theorem~\ref{th:limitDT} for large and small eigenvalues (for discrete time two-colour urns). In~\cite{CMP13}, it is used as in the present paper, to prove the unicity of the solution of a two-equation system, in the study of large two--colour P{\'o}lya urns. In~\cite{Janson04}, it is also used to prove the unicity of the solution of system~\eqref{eq:dcoul_sysCT}: therefore, we will only develop the proof for the discrete case. Similar proofs can be found in~\cite{KN13} or~\cite{CMP13}.

Let $\mathcal{M}_2$ be the space of complex-valued square integrable probability measures. For all $A\in\mathbb{C}$, let $\mathcal{M}^{\mathbb{C}}_2(A)$ be the subspace of measures in $\mathcal{M}_2$ with mean $A$.
We consider the Wasserstein distance as follows: for all $\mu$, $\nu$ two measures in $\mathcal{M}^{\mathbb{C}}_2(A)$,
\[d_W(\mu,\nu) = \inf_{X\sim\mu, Y\sim\nu} \|X-Y\|_2,\]
where $\|\cdot\|_2$ is the $L^2$-norm on $\mathbb{C}$.

For all $A_1,\ldots, A_d\in\mathbb{C}$. 
Let us denote by $\bigtimes_{i=1}^d \mathcal{M}^{\mathbb{C}}_2(A_i)$ 
the Cartesian product of the spaces $\mathcal{M}^{\mathbb{C}}_2(A_i)$.
We define the Wasserstein distance on this space as follows: 
for all $\boldsymbol{\mu}=(\mu_1,\ldots,\mu_d)$ and $\boldsymbol{\nu}=(\nu_1,\ldots,\nu_d)$ 
two elements of $\bigtimes_{i=1}^d \mathcal{M}^{\mathbb{C}}_2(A_i)$,
\[d(\boldsymbol{\mu},\boldsymbol{\nu}) = \max_{1\leq i\leq d}\{d_W(\mu_i,\nu_i)\}.\]
We know that $(\mathcal{M}_2^{\mathbb{C}}(A),d_W)$ and thus $\bigtimes_{i=1}^d \mathcal{M}^{\mathbb{C}}_2(A_i)$ are complete metric spaces (see for example~\cite{Dudley02}). 

The random vector $(W^{DT}_{\boldsymbol{e_1}},\ldots,W^{DT}_{\boldsymbol{e_d}})$ is a solution of System~\eqref{eq:dcoul_sysDT}:
\[W_{\boldsymbol{e_c}} \EgalLoi \sum_{i=1}^d \sum_{k=\gamma_{i-1}^{(c)}+1}^{\gamma_i^{(c)}} \left(V_k^{(c)}\right)^{\nicefrac{\lambda}{S}} W_{\boldsymbol{e_i}}^{(k)}.\]
For all $\boldsymbol{\mu} = (\mu_1,\ldots,\mu_d)\in \bigtimes_{i=1}^d \mathcal{M}^{\mathbb{C}}_2(m_i)$, 
for all $c\in\{1,\ldots,d\}$, let
\[K_c(\boldsymbol{\mu}) = \mathcal{L}\left(\sum_{i=1}^{d} \sum_{k=\gamma_{i-1}^{(c)}+1}^{\gamma_i^{(c)}} 
\left(V_k^{(c)}\right)^{\nicefrac{\lambda}{S}} X_i^{(k)}\right),\]
where $\gamma_0 = 0$, $\gamma_i^{(c)} = \sum_{j\leq i}(\tilde a_{c,j}+\delta_{c,j})$ and for all $i\in\{1,\ldots,d\}$, 
the $(X_i^{(k)})_{1\leq i\leq d}$ are independent random variables, 
independent of each other and of vector $V$, 
which is Dirichlet-distributed of parameter $\left(\pi_1^{(c)}, \ldots \pi_{\gamma_d^{(c)}}^{(c)}\right)$,
and, for all $1\leq i\leq d$ and $1\leq k\leq \gamma_d^{(c)}$,
\begin{equation}\label{eq:pi}
\gamma_{i-1}^{(c)}+1\leq k\leq \gamma_i^{(c)} \Rightarrow X_i^{(k)} \sim \mu_i \text{ and }\pi_k^{(c)}= \theta_i/S.
\end{equation}
We define the function $K$ as
\[K(\boldsymbol{\mu}) = (K_1(\boldsymbol{\mu}),\ldots,K_d(\boldsymbol{\mu})),\]
and prove the following result:
\begin{proposition}
\label{prop:dcoul_contractionDT}
For all large eigenvalue $\lambda$ of the replacement matrix $R$,
for all $\bs{A} = (A_1,\ldots,A_d)\in\mathbb{C}^d$, denote by $\Psi(\bs A)$ the vector whose coordinates are given by
\[\Psi(\bs A)_i = \frac{A_i}{\lambda+\theta_i} \quad \text{ for all }i\in\{1, \ldots, d\},\]
where $\lambda+\theta_i\neq 0$ since $\RRe\lambda>S/2>0$ and $\theta_i\geq 1$.
\begin{enumerate}[(i)]
\item For all $\bs{A} = (A_1,\ldots,A_d)\in\mathbb{C}^d$ such that $\Psi(\bs A)\in\Ker(R-\lambda I_d)$,
the function $K$ is an application from $\bigtimes_{i=1}^d \mathcal{M}^{\mathbb{C}}_2(A_i)$ into itself.

\item Moreover, the law of $(W^{DT}_{\boldsymbol{e_1}},\ldots,W^{DT}_{\boldsymbol{e_d}})$ is the unique square-integrable 
solution of~\eqref{eq:dcoul_sysDT} at fixed mean.
\end{enumerate}
\end{proposition}

\begin{remark}
If we assume $\mathtt{(T_{-1})}$ in addition, then, remark that $\Psi(\bs A)\in \Ker(R-\lambda I_d)$ if and only if $\bs A\in\Ker(R-\lambda I_d)$.
\end{remark}

\begin{proof}
$(i)$
First remark that
\[
\mathbb{E} K_c(\boldsymbol{\mu}):=
\int_{\mathbb C} x \, dK_c(\bs \mu)(x)
= \sum_{i=1}^d \sum_{p=\gamma_{i-1}^{(c)}+1}^{\gamma_i^{(c)}} \mathbb{E} \left(V_p^{(c)}\right)^{\nicefrac{\lambda}{S}} \mathbb{E} X_i^{(k)}
\]
because $(V_1,\ldots,V_{S+1})$ is independent of $(X_i^{(1)},\ldots,X_i^{(S+1)})_{1\leq i\leq d}$. Since, for all $p\in\{1,\ldots, S+1\}$, for all $i\in\{1,\ldots,d\}$, $\mathbb{E}X_i^{(p)} = A_i$, we have
\begin{align*}
\mathbb{E} K_c(\boldsymbol{\mu})
&= \sum_{i=1}^d A_i \sum_{p=\gamma_{i-1}^{(c)}+1}^{\gamma_i^{(c)}} \mathbb{E}\left(V_p^{(c)}\right)^{\nicefrac{\lambda}{S}}\\
&= \sum_{i=1}^d A_i \frac{\gamma_i^{(c)}-\gamma_{i-1}^{(c)}}{1+\lambda\theta_i^{-1}}
= \sum_{i=1}^d A_i \frac{\tilde a_{c,i}+\delta_{c,i}}{1+\lambda\theta_i^{-1}}\\
&= \sum_{i=1}^d A_i \frac{a_{c,i}+\delta_{c,i}\theta_i}{\theta_i+\lambda}
= \sum_{i=1}^d (a_{c,i}+\delta_{c,i}\theta_i) B_i,
\end{align*}
where $\bs B = \Psi(\bs A)$. 
The above calculations are true 
because $\left(V_p^{(c)}\right)^{\nicefrac{\lambda}{S}} \EgalLoi U^{\nicefrac{\lambda}{\theta_c}}$ 
for all $p\in\{1,\ldots, S+1\}$, where $U$ is a random variable uniformly distributed on $[0,1]$,
and because $\RRe\lambda>\nicefrac{S}2$, which implies $\lambda \neq -\theta_c$ 
since $\theta_c\geq 1$ (see Equation~\eqref{eq:theta_i}).
Since~$\lambda$ is an eigenvalue of~$R$, 
and $\boldsymbol{B}= (B_1,\ldots, B_d)\in\Ker(R-\lambda I_d)$, we have
\[\sum_{i=1}^d a_{c,i} B_i= \lambda B_c\]
for all $1\leq c\leq d$. It implies
\[\mathbb{E} K_c(\boldsymbol{\mu}) = (\lambda+\theta_c)B_c = A_c\]
for all $\boldsymbol{\mu}\in \bigtimes_{i=1}^d \mathcal{M}^{\mathbb{C}}_2(A_i)$.
Moreover $K(\boldsymbol{\mu})$ is square-integrable, which implies that $K$ is indeed a function from $\bigtimes_{i=1}^d \mathcal{M}^{\mathbb{C}}_2(A_i)$ into itself, for all $\boldsymbol{A}$ such that $\Psi(\bs A)\in\Ker(R-\lambda I_d)$.

$(ii)$
Let $\bs \mu = (\mu_1,\ldots,\mu_d)\in\bigtimes_{i=1}^d \mathcal{M}^{\mathbb{C}}_2(A_i)$ and 
$\bs \nu = (\nu_1,\ldots,\nu_d)\in\bigtimes_{i=1}^d \mathcal{M}^{\mathbb{C}}_2(A_i)$ 
be two solutions of System~\eqref{eq:dcoul_sysDT}, meaning that
\[K\bs \mu = \bs \mu \quad \text{ and }\quad K\bs \nu = \bs \nu.\]

Let us prove that $d(\bs \mu, \bs \nu) = 0$, using the total variance law: 
it is enough to prove that, for all $i\in\{1, \ldots, d\}$, $d_W(\mu_i, \nu_i) = 0$.
Kantorovitch-Rubinstein's theorem implies that, for all $1\leq i\leq d$, there exists
two random variables $\mathcal X_i\sim \mu_i$ and $\mathcal Y_i\sim\nu_i$ such that 
$d_W(\mu_i, \nu_i) = \|\mathcal X_i - \mathcal Y_i\|_2$.
Fix $1\leq c\leq d$, and let $V$ be a Dirichlet-distributed random variable of parameter 
$\left(\pi_1^{(c)}, \ldots \pi_{\gamma_d^{(c)}}^{(c)}\right)$ (as in Equation~\eqref{eq:pi}),
$X = (X_1^{(c)}, \ldots, X_{\gamma_d^{(c)}}^{(c)})$ and $Y = (Y_1^{(c)}, \ldots, Y_{\gamma_d^{(c)}}^{(c)})$
be two sequences of independent random variables such that $X$ and $Y$ are both independent of $V$ 
(but not necessarily independent of each other), and,
for all $1\leq i\leq d$ and $1\leq k \leq \gamma_d^{(c)}$,
\[\gamma_{i-1}^{(c)}+1\leq k\leq \gamma_i^{(c)} \Rightarrow 
\mathcal L(X_k^{(c)}, Y_k^{(c)}) = \mathcal L(\mathcal X_i, \mathcal Y_i) \text{ and }
\pi_k^{(c)} = \theta_i/S.\]
We thus have, for all $c\in\{1,\ldots, d\}$,
\begin{align*}
d_W(\mu_c,\nu_c)^2
&= d_W(K_c(\boldsymbol{\mu}),K_c(\boldsymbol{\nu}))^2
\leq \left\|\sum_{i=1}^d \sum_{k=\gamma_{i-1}^{(c)}+1}^{\gamma_i^{(c)}} \left(V_k^{(c)}\right)^{\nicefrac{\lambda}{S}} (X_i^{(k)}-Y_i^{(k)})\right\|_2^2&\\
&\leq \sum_{i=1}^d \sum_{k=\gamma_{i-1}^{(c)}+1}^{\gamma_i^{(c)}} 
\mathbb E\left|\left(V_k^{(c)}\right)^{\nicefrac{2\lambda}{S}}\right| \mathbb E\big|(X_i^{(k)}-Y_i^{(k)})\big|^2\\
&= \sum_{i=1}^d \mathbb E|\mathcal X_i-\mathcal Y_i|^2
\sum_{k=\gamma_{i-1}^{(c)}+1}^{\gamma_i^{(c)}} 
\mathbb E\left|\left(V_k^{(c)}\right)^{\nicefrac{2\lambda}{S}}\right|,
\end{align*}
implying that
\[
d_W(K_c(\boldsymbol{\mu}),K_c(\boldsymbol{\nu}))^2
\leq \sum_{i=1}^d \frac{\gamma_i^{(c)}-\gamma_{i-1}^{(c)}}{2\RRe \lambda \theta_i^{-1} +1}\, d_W(\mu_i, \nu_i)^2
= \sum_{i=1}^d  \frac{\tilde a_{c,i}+\delta_{c,i}}{2\RRe\lambda \theta_i^{-1} +1}\,d_W(\mu_i,\nu_i)^2.
\]
If we let $\Delta_i = \frac{d_W(\mu_i, \nu_i)^2}{\theta_i+2\RRe\lambda}$, for all $i\in\{1,\ldots, d\}$, we get that,
for all $1\leq c\leq d$,
\[2\RRe\lambda\Delta_c \leq \sum_{i=1}^d a_{i,c} \Delta_i = (R\Delta)_c.\]
Let $\bs v = (v_1, \ldots, v_d)$ be a horizontal vector with positive entries such that $\bs vR = S\bs v$, then
\[\bs v(R\Delta) = S \bs v\Delta.\]
The existence of such a $\bs v$ is a consequence of the irreducibility of the urn, namely hypothesis $\mathtt{(S)}$, as explained in~\cite{Janson04}.
It implies that,
\[S \sum_{c=1}^d v_c \Delta_c 
=\sum_{c=1}^d v_c (R\Delta)_c 
\geq 2\RRe\lambda \sum_{c=1}^d v_c \Delta_c.\]
Since $\nicefrac{\RRe\lambda}{S}>\nicefrac12$, this last inequality implies that
\[\sum_{c=1}^d v_c\Delta_c = 0,\] and thus, 
by positivity of the $v_c$ and non negativity of the $\Delta_c$, we get,
for all $c\in\{1, \ldots, d\}$, $\Delta_c=0$.
It thus implies that, for all $c\in\{1, \ldots, d\}$,
\[\mu_c=\nu_c, \quad\text{ and thus }\quad \bs \mu=\bs \nu,\]
which concludes the proof.
\end{proof}

\subsection{Decomposition in discrete time}
The argument used to prove Proposition~\ref{thm:dcoul_dislocDT} can also be used to prove the following result from Proposition~\ref{th:dcoul_decompCT} and Theorem~\ref{thm:dcoul_moments}:
\begin{proposition}\label{th:dcoul_decompDT}
Under assumptions $\mathtt{(B)}$, $\mathtt{(T)}$ and $\mathtt{(S)}$,
for all initial composition~$\bs \alpha$,
\[W^{DT}_{\boldsymbol{\alpha}} 
\EgalLoi 
\sum_{c=1}^{d} \sum_{p=\beta_{c-1}+1}^{\beta_c} 
Z_p^{\nicefrac{\lambda}{S}} W_{\boldsymbol{e_c}}^{(p)},\]
where $\beta_0 = 0$, $\beta_i = \sum_{j\leq i} \tilde \alpha_j$;
where the $W_{\boldsymbol{e_c}}^{(p)}$ are independent copies of $W^{DT}_{\boldsymbol{e_c}}$, independent of each other;
and where $Z = (Z_1,\ldots,Z_{\beta_d})$ is a Dirichlet-distributed random vector independent of the $W$ and of parameter~$\bs \eta$ given by
\[\bs \eta_k = \nicefrac{\theta_i}{S}\quad \text{ if } \beta_{i-1}<k\leq \beta_i.\]
\end{proposition}

\begin{remark}
If we assume $\mathtt{(T_{-1})}$ instead of $\mathtt{(T)}$ in the theorem above, we obtain the following equation:
for all initial composition $\bs \alpha$,
\[W^{DT}_{\boldsymbol{\alpha}} 
\EgalLoi \sum_{c=1}^{d} \sum_{p=\beta_{c-1}+1}^{\beta_c} Z_p^{\nicefrac{\lambda}{S}} W_{\boldsymbol{e_c}}^{(p)},\]
where $\beta_0 = 0$, $\beta_i = \sum_{j\leq i}\alpha_j$;
where $Z = (Z_1,\ldots,Z_{\beta_d})$ is a Dirichlet-distributed random vector of parameter~$\left(\frac1{S},\ldots,\frac1{S}\right)$;
and where the $W_{\boldsymbol{e_c}}^{(p)}$ are independent copies of $W^{DT}_{\boldsymbol{e_c}}$, independent of each other and of $Z$.
\end{remark}

\subsection{Tree structure in discrete time}\label{sub:alt_proof}
Propositions~\ref{th:dcoul_decompDT} and~\ref{thm:dcoul_dislocDT} can be proven 
from scratch by an analogue of the proof of Propositions~\ref{th:dcoul_decompCT} and~\ref{thm:dcoul_dislocCT}, 
using the underlying tree structure of the urn.
The analysis of the tree structure in discrete-time is however more intricate since the 
subtrees of the considered forest are not independent.
Such a tree decomposition in discrete-time is already proposed in~\cite{CMP13} (for two-colour urns) or~\cite{KN13} but both papers assume $\mathtt{(B)}$, $\mathtt{(I)}$ and $\mathtt{(T_{-1})}$.
We will develop this alternative proof of Proposition~\ref{th:dcoul_decompDT} 
because it happens to be more complicated due to the possibly negative diagonal coefficient of the replacement matrix.

\begin{proof}[Alternative proof of Proposition~\ref{th:dcoul_decompDT}]
The discrete-time urn process can be seen as a forest whose leaves can be of $d$ different colours: 
we refer to Figure~\ref{fig:decomp_DT}\footnote{Figure~\ref{fig:decomp_DT} is an example of a two-colour urn. 
The example taken is chosen for it simplicity even if the chosen urn only has small Jordan blocks: 
it is called a \emph{small} urn in the literature. 
It does not affect the arguments developed in the proof.}. 
A time zero, the forest is composed of $\tilde \alpha_i$ roots of colour $i$ 
(for all $i\in\{1,\ldots, d\}$), 
each of these roots contains $\theta_i$ balls of colour $i$. 
At each step, we pick up uniformly at random a ball: 
this ball belongs to a leaf. 
The picked leaf then becomes an internal node, 
which has $\gamma_d^{(c)}$ children, 
amongst them $\tilde a_{i,j} + \delta_{i,j}$ contain $\theta_j$ balls of colour $j$ 
(for all $j\in\{1,\ldots, d\}$) if the picked up leaf was of colour $i$.

\begin{figure}
\caption[bla]{A realisation at time $n=7$ of the forest associated to 
the urn process with initial composition $^t\!(2,2)$ 
and with replacement matrix 
$R=\begin{pmatrix} -2 & 4 \\ 2 & 0\end{pmatrix}$.
Remark that $S_1(7) = 3$, $S_2(7) = 3$ and $S_3(7) = 1$.
Moreover, $D_1(7) = 8$, $D_2(7) = 3$ and $D_3(7) = 1$.
}
\label{fig:decomp_DT}
\begin{center}
\includegraphics[width=.8\textwidth]{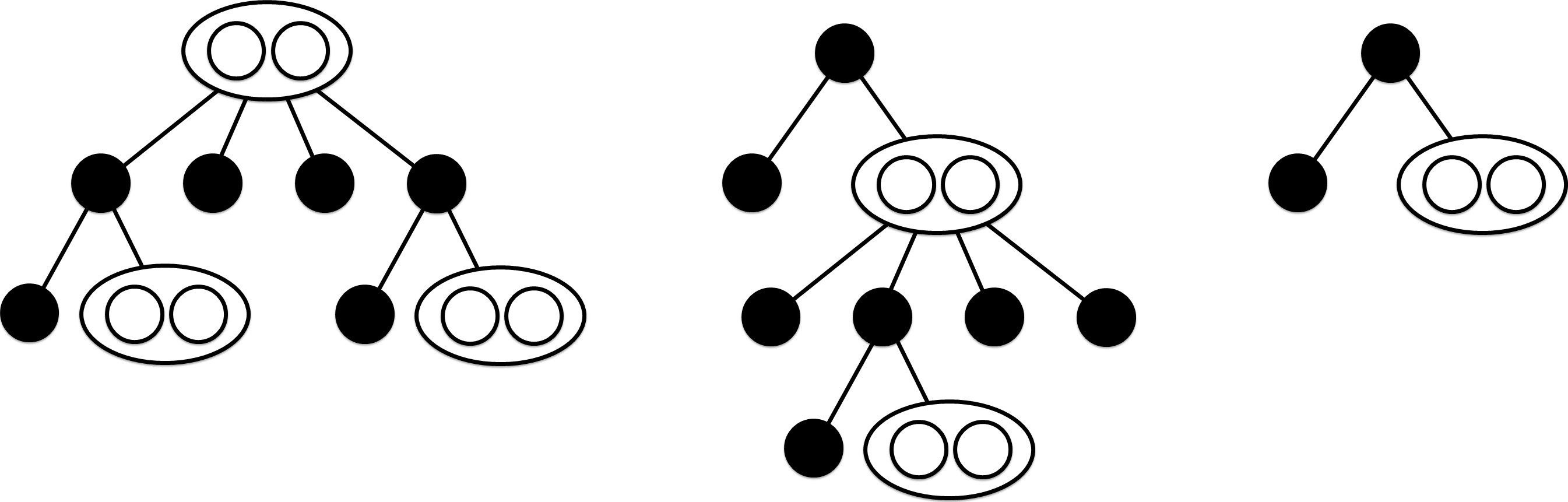}
\end{center}
\end{figure}

The composition of the urn is thus described by the set of leaves of the forest.
Let us number the subtrees of the forest from trees rooted by colour $1$ up to trees rooted by colour $d$.
If we denote by $D_p(n)$ the number of balls in leaves of the $p^{\text{th}}$ subtree of the forest, 
then this $p^{\text{th}}$ subtree of the forest at time $n$ represents the composition vector 
of an urn process with inital composition of cardinal $\theta_c$ if $\beta_{c-1}<p\leq \beta_c$,
taken at \emph{internal time} $S_p(n) = \frac{D_p(n)-\theta_c}{S}$.
Indeed, the \emph{internal time} in the $p^{\text{th}}$ tree is the number of its internal nodes; 
the fact that the urn is balanced means that all internal nodes of the tree have given birth to $S+1$ balls,
which give the above relationship between leaves and internal nodes in the $p^{\text{th}}$ subtree. 

We thus have
\begin{equation}
\label{eq:decompDT}
U_{\boldsymbol{\alpha}}(n) \EgalLoi 
\sum_{c=1}^{d} \sum_{p=\beta_{c-1}+1}^{\beta_c} 
U_{\boldsymbol{e_c}}^{(p)}\left(\frac{D_p(n)-\omega_p}{S}\right),
\end{equation}
where $\beta_0 = 0$, where for all $c\geq 1$, $\beta_c = \sum_{i=1}^{c}\tilde \alpha_i$;
where $\omega_p = \theta_i$ if $\beta_{i-1}<p\leq \beta_i$;
and where the urn processes $U_{\boldsymbol{e_c}}^{(p)}$ are independent copies of the process $U_{\boldsymbol{e_c}}$, 
independent of each other.
We are thus interested in the asymptotic behaviour of the vector $(D_1(n),\ldots,D_{\beta_d}(n))$ when $n$ grows to infinity. 
This vector happens to be the composition vector of a $\left(\beta_d = \sum_{i=1}^d \tilde \alpha_i\right)$--colour P{\'o}lya urn with initial composition $\bs\omega = {}^t \!(\omega_1, \ldots, \omega_{\beta_d})$,
where
\[\omega_k = \theta_i \quad \text{ if } \beta_{i-1}<k\leq \beta_i.\]
Indeed, forget the initial colouring of the leaves and colour the leaves of the $i$th subtree with colour $i$; 
at each step, a leaf of the forest picked up uniformly at random becomes an internal node and gives birth to $S+1$ leaves of its same colour. 
Such a diagonal urn is called a P{\'o}lya-Eggenberger urn and has been long studied in the literature. 
We can for example cite this result by \cite{Athreya69} 
(for a complete proof, see \cite{Bertoin06} or \cite{CMP13}): 

\begin{theorem}\label{thm:Dirichlet}
Let $p$ and $K$ be two positive integers and $(D_1(n),\ldots,D_p(n))$ the composition vector at time $n$ of an urn process of initial composition $\bs\nu \!= \!{}^t \!(\nu_1, \ldots, \nu_p)$ and with replacement matrix $K I_{p}$, then, 
asymptotically when $n$ tends to infinity, almost surely,
\[\frac{1}{nK}\left(D_1(n),\ldots,D_{p}(n)\right)\to \bs Z = (Z_1,\ldots,Z_{p})\]
where $\bs Z$ is a Dirichlet-distributed random vector of parameter 
$\left(\frac{\nu_1}{K},\ldots,\frac{\nu_p}{K}\right)$.
\end{theorem}

Thus, projecting Equation~\eqref{eq:decompDT} onto $E$ via $\pi_E$, renormalising it by $n^{\nicefrac{\lambda}{S}} \ln^{\nu} n$ and taking the limit when $n\to+\infty$ gives Proposition~\ref{th:dcoul_decompDT}.
\end{proof}

\section{Densities -- Proof of Theorem~\ref{th:density}}\label{sec:d_coul_density}
In this section, we prove Theorem~\ref{th:density} via an analysis of the Fourier transforms of the random variables $W^{CT}$ and $W^{DT}$. We generalise the method developed by \cite{Liu98,Liu01} for smoothing equations with positive solutions and refer to~\cite{CLP13} for a similar proof in the case of $m$-ary trees where a complex fixed point equation (but not a system) is studied. 
A similar result is also proved for the two--colour case in~\cite{CMP13}, 
but due to both the higher dimension and the weaker tenability condition,
the present proof follows a different route.

The strategy of the proof is the following. First focus on the discrete time case, i.e. on $W^{DT}$:
\begin{itemize}
\item We prove that the support of $W^{DT}_{\bs {e_i}}$ contains some non-discrete set.
\item It implies that the Fourier transform of $W^{DT}$ is integrable and thus invertible (see Lemmas~\ref{lem:psi<1},~\ref{lem:Riemann} and~\ref{lem:rho}), so that the variable $W^{DT}$ has a density on $\mathbb{C}$.
\item We finally deduce from the existence of a density and from Lemma~\ref{lem:support} that the support of $W^{DT}$ is the whole complex plane.
\end{itemize}
Via the martingale connection~\eqref{eq:mcCT}, one can consequently infer that, 
for all $c\in\{1,\ldots,d\}$ $W^{CT}_{\boldsymbol{e_c}}$ has a density on $\mathbb{C}$. 
Propositions~\ref{th:dcoul_decompDT} and~\ref{th:dcoul_decompCT} permit to generalise to any initial composition.

We assume $\mathtt{(B)}$, $\mathtt{(T)}$ and $\mathtt{(I)}$ and recall that the eigenvalue $\lambda$, 
the Jordan subspace $E$, and its dimension $\nu+1$ are fixed.
\begin{lemma}\label{lem:esperance}
There exists $1\leq c\leq d$ such that $\mathbb{E}W_{\bs{e_c}}^{DT}\neq 0$.
\end{lemma}
\begin{proof}
It is known (see for example~\cite{Pouyanne08}) that, for all $c\in\{1, \ldots, d\}$,
\begin{align*}
&\frac1{\nu!}W^{DT}_{\bs{e_c}}\bs v \\
&= \lim_{n\to+\infty} \left(1+{{}^t \!R_{|E}}\right)^{-1} \left(1+\frac{{}^t \!R_{|E}}{1+S}\right)^{-1} \cdots \left(1+\frac{{}^t \!R_{|E}}{1+(n-1)S}\right)^{-1} \pi_E(U_{\bs{e_c}}(n)),
\end{align*}
as a martingale almost sure limit, where $^t\!R_{|E}$ stands for the restriction of 
the endomorphism induced by~$R$ to the Jordan subspace~$E$ 
(note that the unique eigenvalue of $R$ is thus $\lambda$ and since $\RRe\lambda>S/2$, the inverses above are well defined).
Thus, $\mathbb{E}W^{DT}_{\bs{e_c}}\bs v = \nu! \pi_E(\bs{e_c})$.
Note that, if $\mathbb{E}W^{DT}_{\bs{e_c}}\bs v = \nu! \pi_E(\bs{e_c}) = 0$, for all $1\leq c\leq d$,
then, $\pi_E =0$, which is impossible.
\end{proof}

\begin{lemma}\label{lem:support}
For all $\boldsymbol{z} = (z_1, \ldots, z_d)\in\bigtimes_{i=1}^d \supp(W^{DT}_{\boldsymbol{e_i}})$, 
for all $c\in\{1, \ldots, d\}$,
for all $(v_1,\ldots,v_{\gamma_d^{(c)}})$ in the support of a Dirichlet-
distributed random vector of parameter $\bs \omega= (\omega_1, \ldots, \omega_{\gamma_d^{(c)}})$,
where $\omega_i = \nicefrac{\theta_i}{S}$ for all $i\in\{1, \ldots, c\}$, we have
\[\sum_{i=1}^d \sum_{k=\gamma_{i-1}^{(c)}+1}^{\gamma_i^{(c)}} v_k^{\nicefrac{\lambda}{S}} z_i \in \Supp(W^{DT}_{\bs{e_c}}).\]
\end{lemma}
Note that the support of the Dirichlet distribution of parameter $\bs \omega$ is
\[\Big\{(v_1, \ldots, v_{\gamma_d^{(c)}})\in [0,1]^{\gamma_d^{(c)}} \colon \sum_{k=1}^{\gamma_d^{(c)}}v_k = 1\Big\}.\]
\begin{proof}
Recall that for a given complex random variable $Z$, for all $z\in\mathbb{C}$,
\[z\in \supp(Z) \Leftrightarrow \forall \varepsilon>0, \mathbb{P}(|Z-z|<\varepsilon)>0.\]
Let $\boldsymbol{z}$, $c$ and $(v_1,\ldots,v_{\gamma_d^{(c)}})$ satisfying the hypothesis of the lemma and let $\eta>0$. We have
\begin{align*}
\left|W^{DT}_{\boldsymbol{e_c}} - \sum_{i=1}^d \sum_{k=\gamma_{i-1}^{(c)}+1}^{\gamma_i^{(c)}} v_k^{\nicefrac{\lambda}{S}} z_i\right|
\EgalLoi \left|\sum_{i=1}^d \sum_{k=\gamma_{i-1}^{(c)}+1}^{\gamma_i^{(c)}} \!\!
\left(V_k^{(c)}\right)^{\!\nicefrac{\lambda}{S}} W_{\boldsymbol{e_i}}^{(k)}
- \sum_{i=1}^d\sum_{k=\gamma_{i-1}^{(c)}+1}^{\gamma_i^{(c)}} \!\!
v_k^{\nicefrac{\lambda}{S}} z_i\right|&\\
\EgalLoi \left|\sum_{i=1}^d \sum_{k=\gamma_{i-1}^{(c)}+1}^{\gamma_i^{(c)}}\!\!
 \left(V_k^{(c)}\right)^{\nicefrac{\lambda}{S}} (W_{\boldsymbol{e_i}}^{(k)}-z_i)
- \sum_{i=1}^d\sum_{k=\gamma_{i-1}^{(c)}+1}^{\gamma_i^{(c)}} \!\!
\left(v_k^{\nicefrac{\lambda}{S}}-\left(V_k^{(c)}\right)^{\nicefrac{\lambda}{S}}\right) z_i\right|&,
\end{align*}
where $V = (V_1^{(c)}, \ldots, V_{\gamma_d^{(c)}}^{(c)})$ and the $W_{\bs e_i}^{(k)}$ are defined as in Theorem~\ref{thm:dcoul_dislocDT}.
Thus, since $(z\mapsto z^{\nicefrac{\lambda}{S}})$ is continuous on $\mathbb{C}$, we have, with positive probability,
\[|W_{\bs e_i}^{(k)}-z_i|< \eta \quad \text{ and }\quad
\big|v_k^{\nicefrac{\lambda}{S}}-\big(V_k^{(c)}\big)^{\nicefrac{\lambda}{S}}\big| < \eta,\]
which imply
\begin{align*}
\left|W^{DT}_{\boldsymbol{e_c}} - \sum_{i=1}^d \sum_{p=\gamma_{i-1}^{(c)}+1}^{\gamma_i^{(c)}} v_p^{\nicefrac{\lambda}{S}} z_i\right|
&\leq \eta \sum_{p=1}^{\gamma_d^{(c)}} \left|\left(V_p^{(c)}\right)^{\nicefrac{\lambda}{S}}\right| + \eta \sum_{i=1}^d |z_i|\\
&\leq (\gamma_d^{(c)}+\|z\|_1)\eta \leq (S+1+\|z\|_1)\eta,    
\end{align*}
where $\|z\|_1=\sum_{i=1}^d |z_i|$. For all $\varepsilon>0$, we fix $\eta = \frac{\varepsilon}{S+1+\|z\|_1}$ to conclude the proof.
\end{proof}

\begin{lemma}\label{lem:densite}
There exists a non zero $z_0$ such that, for all $t\in(0,1)$,
\[(t^{\nicefrac{\lambda}{S}}+(1-t)^{\nicefrac{\lambda}{S}})z_0\in \bigcap_{i=1}^d \Supp(W_{\bs{e_i}}^{DT}).\]
\end{lemma}

\begin{proof} 
Thanks to Lemma~\ref{lem:esperance}, there exists a colour $1\leq c\leq d$ 
such that there exists $z_0\neq 0\in\Supp(W_{\bs{e_c}}^{DT})$. 
Applying Lemma~\ref{lem:support}, for all $i\in\{1, \ldots, d\}$ 
such that $a_{c,i}\neq 0$, $z_0\in\Supp(W_{\bs{e_i}}^{DT})$.
Iterating this argument allows us to conclude that, 
for all $i\in\{1, \ldots, d\}$ dominated by $c$, $z_0\in\Supp(W_{\bs{e_i}}^{DT})$.
By Assumption $\mathtt{(I)}$, all colours are dominating, and thus $c$ is dominating, meaning that
\[z_0\in \bigcap_{i=1}^d \Supp(W_{\bs{e_i}}^{DT}).\]
Therefore, still by applying Lemma~\ref{lem:support}, for all $t\in[0,1]$,
\[(t^{\nicefrac{\lambda}{S}}+(1-t)^{\nicefrac{\lambda}{S}})z_0\in \bigcap_{i=1}^d \Supp(W_{\bs{e_i}}^{DT}).\qedhere\]
\end{proof}

The three following lemmas are proven via very similar arguments as the ones developed in~\cite{CMP13}. 
There is no additional idea to the proof here, except being careful to the slight changes induced by the higher $d\geq 3$ 
and by the weaker assumption $\mathtt{(T)}$ instead of $\mathtt{(T_{-1})}$.
For all $c\in\{1,\ldots,d\}$, let $\phi_c(t) = \mathbb{E} \mathtt{e}^{i\langle t, W^{DT}_{\boldsymbol{e_c}}\rangle}$ 
for all $t\in\mathbb{C}$ and $\psi_c(r) = \sup_{|t|=r} |\phi_c(t)|$.

\begin{lemma}
\label{lem:psi<1}
For all $c\in\{1,\ldots,d\}$, for all $r>0$, $\psi_c(r) <1$.
\end{lemma}

\begin{proof}
We know that $\psi_c(0) = 1$ and that $\psi_c(r)\leq 1$ for all $r\geq 0$. 
Let us assume that there exists $r_c>0$ such that $\psi_c(r_c) = 1$. 
Then, there exists $z_c\in\mathbb{C}$ and $\theta_c\in\mathbb{R}$ such that $|z_c|=r_c$ and 
\[\mathbb{E}\mathtt{e}^{i\langle z_c,W^{DT}_{\bs{e_c}}\rangle} = \mathtt{e}^{i\theta_c}.\]
Thus, the complex random variable $\mathtt{e}^{i\langle z_c,W^{DT}_{\bs{e_c}}\rangle - i\theta_c}$ has mean $1$ and takes its values in the unit disc. It is therefore almost surely equal to $1$, implying that almost surely, $\langle z_c,W^{DT}_{\bs{e_c}}\rangle \in \theta_c + 2\pi\mathbb Z$.
Therefore, there exists $\varpi\in\mathbb Z$ such that, for all $t\in[0,1]$,
\[\langle z_c, (t^{\nicefrac{\lambda}{S}} + (1-t)^{\nicefrac{\lambda}{S}})z_0\rangle
= \langle t^{\nicefrac{\lambda}{S}} + (1-t)^{\nicefrac{\lambda}{S}}, z_c \bar z_0\rangle
= \theta_c+2\pi\varpi.\]
Note that if $w\in\mathbb C\setminus\{0\}$ and $x\in\mathbb R$, then the set $\{z\in \mathbb C \colon \langle z,w\rangle=x\}$ is the line $\nicefrac x w + iw\mathbb R$. Therefore, for all $t\in [0,1]$,
\begin{equation}\label{eq:spiral}
t^{\nicefrac{\lambda}{S}} + (1-t)^{\nicefrac{\lambda}{S}} \in \frac{\theta_c+2\varpi \pi}{\bar z_0 z_c} + i\bar z_0 z_c \mathbb R.
\end{equation}
If $\lambda\in\mathbb R$, then $t^{\nicefrac{\lambda}{S}} + (1-t)^{\nicefrac{\lambda}{S}} \in \mathbb R$ for all $t\in[0,1]$ which implies that both $(\theta_c+2\varpi \pi)/(\bar z_0 z_c)$ and $i\bar z_0 z_c$ are real, which is impossible.
If $\lambda\in\mathbb C\setminus{\mathbb R}$, then, the left-hand side of~\eqref{eq:spiral} is a spiral and thus contains at least three non-aligned points whereas the right-hand side is a line, which is impossible, and thus concludes the proof.
\end{proof}

\begin{lemma}
\label{lem:Riemann}
For all $c\in\{1,\ldots,d\}$, $\lim_{r\to\infty} \psi_c(r) = 0$.
\end{lemma}

\begin{proof}
In view of Equation~\eqref{eq:dcoul_sysDT}, for all $c\in\{1,\ldots,d\}$,
\begin{align*}
\phi_c(t) 
&= \mathbb{E}\exp\left(i\Big\langle t, \sum_{i=1}^d\sum_{k=\gamma_{i-1}^{(c)}+1}^{\gamma_i^{(c)}} 
\left(V_k^{(c)}\right)^{\nicefrac{\lambda}{S}} W_{\bs{e_i}}^{(k)}\Big\rangle\right)\\
&= \mathbb{E}\bigg[\mathbb{E}\bigg[\exp\Big(i\Big\langle t, \sum_{i=1}^d\sum_{k=\gamma_{i-1}^{(c)}+1}^{\gamma_i^{(c)}} \left(V_k^{(c)}\right)^{\nicefrac{\lambda}{S}} W_{\bs{e_i}}^{(k)}\Big\rangle\Big) \bigg | V_1,\ldots, V_{\gamma_d^{(c)}+1}\bigg]\bigg]\\
&= \mathbb{E}\prod_{i=1}^d\prod_{k=\gamma_{i-1}^{(c)}+1}^{\gamma_i^{(c)}} 
\mathbb{E}\left[\mathtt{e}^{i\langle t, \left(V_k^{(c)}\right)^{\nicefrac{\lambda}{S}}W_{\bs{e_i}}^{(k)}\rangle} \bigg | V_1, \ldots, V_{\gamma_d^{(c)}+1}\right],
\end{align*}
implying that
\[\psi_c(r) \leq \mathbb{E}\prod_{i=1}^d \prod_{k={\gamma_{i-1}^{(c)}}+1}^{\gamma_i^{(c)}} 
\psi_i\Big(\big(V_k^{(c)}\big)^{\nicefrac{\RRe\lambda}{S}}r\Big).\]
By Fatou's Lemma, since $\mathbb P(V_k^{(c)}=0) = 0$ for all $1\leq k\leq \gamma_d^{(c)}$, it implies
\[\limsup_{r\to+\infty} \psi_c(r) 
\leq \mathbb{E}\prod_{i=1}^d\big(\limsup_{r\to+\infty}\psi_i(|V_i^{\nicefrac{\lambda}{S}}| r)\big)^{\tilde a_{c,i}+\delta_{c,i}}.
\]
Let us assume that there exists a colour~$c$ such that $\limsup_{r\to\infty}\psi_c(r)=1$.
Let us define $\psi(r) = \max_{1\leq i\leq d}\psi_i(r)$, then $\limsup_{r\to\infty} \psi(r) = 1$.
According to Lemma~\ref{lem:psi<1}, for all $i\in\{1,\ldots,d\}$, $\psi_i(1)<1$, 
and thus, $\psi(1)<1$. 
Let $\varepsilon\in(0,1-\psi(1))$, and define:
\[
\begin{array}{l}
r_1(\varepsilon) = \max\{r\in(0,1), \psi(r) = 1-\varepsilon\}\\
r_2(\varepsilon) = \min\{r>1, \psi(r) = 1-\varepsilon\}.
\end{array}
\]
These definitions are legal since $\psi$ is continuous, $\psi(0) = 1$ and $\limsup_{r\to\infty} \psi(r) = 1$. 
Moreover, we have $\psi(r_1(\varepsilon))=\psi(r_2(\varepsilon)) = 1-\varepsilon$, 
and for all $r\in[r_1(\varepsilon),r_2(\varepsilon)]$, $\psi(r)\leq 1-\varepsilon$. 

We know that for all $c\in\{1,\ldots,d\}$
\[\psi_c(r) \leq \mathbb{E}\prod_{i=1}^d \prod_{k={\gamma_{i-1}^{(c)}}+1}^{\gamma_i^{(c)}} 
\psi_i\Big(\big(V_k^{(c)}\big)^{\nicefrac{\RRe\lambda}{S}}r\Big).\]
In particular, for all $r\geq 0$,
\[\psi_c(r)\leq \mathbb E\psi_c\left(|V_c^{\nicefrac{\lambda}{S}}|r\right),\]
where $V_c$ is Beta-distributed with parameter $(\theta_c/S, 1)$.
For all $1\leq i\leq d$, let $(A^{(i)}_k)_{k\geq 1}$ be a sequence of i.i.d.\ 
random variables having the same law as $|V_i^{\nicefrac{\lambda}{S}}|$. 
Iterating the last identity, we get
\[\psi_c(r) \leq \mathbb{E}\psi_c(r A^{(c)}_1\ldots A^{(c)}_n).\]
Let us define, for all $i\in\{1, \ldots, d\}$, for all integers~$n$
\[\lambda^{(i)}_n(r,\varepsilon) 
= \mathbb{P}(r_1(\varepsilon)\leq r A^{(i)}_1\ldots A^{(i)}_n\leq r_2(\varepsilon)).\]
For all $r\geq 0$, for all $i\in\{1, \ldots, d\}$,
\[\psi_i(r)\leq \mathbb{E}\psi_c(r A^{(i)}_1\ldots A^{(i)}_n)
\leq (1-\varepsilon)\lambda^{(i)}_n(r,\varepsilon)
+ 1- \lambda^{(i)}_n(r,\varepsilon)
\leq 1-\varepsilon\lambda^{(i)}_n(r,\varepsilon).\]
For all $c\in\{1, \ldots, d\}$, we have,
\begin{align*}
\psi_c(r_2(\varepsilon))
& \leq \mathbb{E}\prod_{i=1}^d \prod_{k={\gamma_{i-1}^{(c)}}+1}^{\gamma_i^{(c)}} 
\psi_i\Big(\big(V_k^{(c)}\big)^{\nicefrac{\RRe\lambda}{S}}r_2(\varepsilon)\Big)\\
&\leq \mathbb{E}\prod_{i=1}^d \prod_{k={\gamma_{i-1}^{(c)}}+1}^{\gamma_i^{(c)}}
\Big(1-\varepsilon\lambda_n^{(i)}\big(\big(V_k^{(c)}\big)^{\nicefrac{\RRe\lambda}{S}}r_2(\varepsilon),\varepsilon\big)\Big).
\end{align*}
%where $A^{(i)}\EgalLoi |V_i^{\nicefrac{\lambda}{S}}| = V_i^{\nicefrac{\RRe\lambda}{S}}$.
In view of Lemma~\ref{lem:psi<1}, $r_1(\varepsilon)\to 0$ when $\varepsilon\to 0$, implying that, 
for all $i\in\{1, \ldots, d\}$, when $\varepsilon$ goes to zero,
\[\lambda_n^{(i)}\big(\big(V_k^{(c)}\big)^{\nicefrac{\RRe\lambda}{S}}r_2(\varepsilon),\varepsilon\big)\to \mathbb{P}(0\leq \big(V_k^{(c)}\big)^{\nicefrac{\RRe\lambda}{S}}A^{(i)}_1\ldots A^{(i)}_n \leq 1) = 1.\]
Thus, for all $c\in\{1, \ldots, d\}$,
\[\frac{1-\mathbb{E}\prod_{i=1}^d \prod_{k={\gamma_{i-1}^{(c)}}+1}^{\gamma_i^{(c)}}
\Big(1-\varepsilon\lambda_n^{(i)}\big(\big(V_k^{(c)}\big)^{\nicefrac{\RRe\lambda}{S}}r_2(\varepsilon),\varepsilon\big)\Big) }{\varepsilon}
\to \sum_{i=1}^d (\tilde a_{c,i}+\delta_{c,i})=\gamma_c^{(d)}.\]
Now note that
\begin{align*}
1-\varepsilon 
&= \psi(r_2(\varepsilon)) = \max_{c\in\{1, \ldots, d\}}\psi_c(r_2(\varepsilon))\\
&\leq \max_{c\in\{1, \ldots, d\}} \mathbb{E}\prod_{i=1}^d \prod_{k={\gamma_{i-1}^{(c)}}+1}^{\gamma_i^{(c)}}
\Big(1-\varepsilon\lambda_n^{(i)}\big(\big(V_k^{(c)}\big)^{\nicefrac{\RRe\lambda}{S}}r_2(\varepsilon),\varepsilon\big)\Big) ,
\end{align*}
which implies that
$\min_{c\in\{1, \ldots, d\}}\gamma_d^{(c)} \leq 1$.
Recall that, in view of Equation~\eqref{eq:theta_i},
\[\gamma_d^{(c)} = 1 + \sum_{i=1}^d \tilde a_{c,i} = 1 + \sum_{i=1}^d \frac{a_{c,i}}{\theta_i} > 1,\]
since at least one of the $a_{c,i}$'s is non-zero as their sum is equal to $S\geq 1$.
We have thus reached a contradiction, implying that for all $1\leq c\leq d$,
\[\limsup_{r\to\infty} \psi_c(r) = 1.\qedhere\] 
\end{proof}

\begin{lemma}\label{lem:rho}
For all $c\in\{1, \ldots, d\}$, for all $\rho\in(0,\nicefrac{\theta_c}{\RRe\lambda})$, 
asymptotically when $|t|$ tends to $+\infty$, $\phi_c(t) = \mathcal{O}(|t|^{-\rho})$.
\end{lemma}

\begin{proof}
Let $\varepsilon>0$. 
In view of Lemma~\ref{lem:Riemann}, 
there exists $T>0$ such that, for all $r\geq T$, 
for all $i\in\{1, \ldots, d\}$, $\psi_i(r)\leq \varepsilon$. 
We already have proved that
\[\psi_c(r) \leq \mathbb{E}\prod_{i=1}^d \prod_{k={\gamma_{i-1}^{(c)}}+1}^{\gamma_i^{(c)}} 
\psi_i\Big(\big(V_k^{(c)}\big)^{\nicefrac{\RRe\lambda}{S}}r\Big).\]
Thus, for all $r\geq T$,
\[\psi_c(r) 
\leq \varepsilon \mathbb{E} \psi_c((V_{j(c)}^{(c)})^{\nicefrac{\RRe\lambda}{S}}r) 
+ \mathbb{P}\left(\left|\left(V_{\gamma_c^{(c)}}^{(c)}\right)^{\!\!\nicefrac{\lambda}{S}}r\right|< T\right),\]
where $j(c)$ is any integer from $\{1, \ldots, d\}\setminus\{\gamma_c^{(c)}\}$.
Recall that $V^{(c)}_{\gamma_c^{(c)}}\EgalLoi U^{\nicefrac{S}{\theta_c}}$ and $V^{(c)}_{j(c)}\EgalLoi U^{\nicefrac{S}{\theta_{j(c)}}}$
where $U$ is a uniform random variable on $(0,1)$.
Then, for all $r > T$,
\[\psi_c(r) 
\leq \varepsilon \mathbb{E}\psi_c\big(U^{{\RRe\lambda}/{\theta_{j(c)}}}r\big) 
+ \left(\frac{T}{r}\right)^{\nicefrac{\theta_c}{\RRe\lambda}}.
\]
Said differently, there exists a positive constant $C$ such that, for all $r>T$,
\[\psi_c(r) \leq \varepsilon \mathbb{E}\psi_c\big(U^{{\RRe\lambda}/{\theta_{j(c)}}} r\big)
+ C \left(\frac{1}{r}\right)^{\rho},\]
for all $\rho\in(0,\nicefrac{\theta_c}{\RRe\lambda})$.
For such $\rho$ (and actually for any $\rho$), 
$\mathbb{E}U^{-{\rho\RRe\lambda}/{\theta_{j(c)}}}<+\infty$, 
and we can thus apply a Gronwall-type lemma (\cite[Lemma~4.1]{Liu99}), which implies
\[\psi_c(r)\leq \frac{Cr^{-\rho}}{1-\varepsilon \mathbb{E}U^{-{\rho\RRe\lambda}{\theta_{j(c)}}}},\]
provided that $\varepsilon$ satisfies $1-\varepsilon \mathbb{E}U^{-{\rho\RRe\lambda}/{\theta_{j(c)}}}>0$.
\end{proof}

\begin{proposition}\label{prop:density}
If $\lambda\in\mathbb{C}\setminus\mathbb{R}$, the distribution of $W^{DT}_{\bs{e_c}}$ admits a density on $\mathbb{C}$ and its support is $\mathbb{C}$, for all $c\in\{1, \ldots, d\}$.

If $\lambda\in\mathbb{R}$, the distribution of $W_{\bs{e_c}}^{DT}$ admits a density on $\mathbb{R}$ and its support is $\mathbb{R}$, for all $c\in\{1, \ldots, d\}$.
\end{proposition}

\begin{proof}
Let us apply arguments already used in~\cite[page~22]{CMP13} to prove that $\phi_c$ is integrable for all $c\in\{1, \ldots, d\}$. Recall that, for all $r\geq 0$,
\[\psi_c(r) \leq \mathbb{E}\prod_{i=1}^d \prod_{k={\gamma_{i-1}^{(c)}}+1}^{\gamma_i^{(c)}} 
\psi_i\Big(\big(V_k^{(c)}\big)^{\nicefrac{\RRe\lambda}{S}}r\Big).\]
In view of Lemma~\ref{lem:rho},  there exists a constant $\kappa>0$ such that, 
for all $i\in\{1, \ldots, d\}$, for all $\rho_i\in(0, \nicefrac{\theta_i}{\RRe \lambda})$, for all $|t|$ large enough,
$|\phi_i(t)|\leq \kappa |t|^{-\rho_i}$.
Let $\eta_c = \sum_{i=1}^d (\tilde a_{c,i}+\delta_{c,i})\rho_i$, 
for all $r$ large enough, we have
\[\psi_c(r)\leq \frac{\kappa^{\gamma_d^{(c)}}}{r^{\eta_c}} 
\mathbb{E}\prod_{i=1}^d \prod_{k={\gamma_{i-1}^{(c)}}+1}^{\gamma_i^{(c)}} 
\big(V_k^{(c)}\big)^{-\nicefrac{\rho_i\RRe\lambda}{S}}
= \mathcal O\left(r^{-\eta_c}\right).\]
Note that, for all $1\leq c\leq d$, in view of Equation~\eqref{eq:theta_i},
\[\sum_{i=1}^d (\tilde a_{c,i}+\delta_{c,i}) \frac{\theta_i}{S} 
= \sum_{i=1}^d \frac{a_{c,i}+ \theta_i\delta_{c,i}}{S}
= 1+\frac{\theta_c}{S} >1,\]
implying that we can choose $\rho_1, \ldots, \rho_d$ such that $\rho_i < \theta_i/S$ and $\eta_c > 1$.
Therefore, $\phi_c$ is integrable and thus that $W_{\bs{e_c}}^{DT}$ admits a 
bounded and continuous density on $\mathbb{C}$, for all $c\in\{1, \ldots, d\}$.

Fix a dominating colour $c$ (actually, by Assumption $\mathtt{(I)}$, 
all colours are dominating, but for this argument, $\mathtt{(S)}$ would be enough).
Since $W_{\bs e_c}^{DT}$ admits a continuous density, we know that there exists $z_0\in \mathbb C$ and $\varepsilon>0$ 
such that the opened ball centred in $z_0$ of radius $\varepsilon$, which we denote by $\mathcal B(z_0, \varepsilon)$,
is contained in $\Supp (W_{\bs{e_c}}^{DT})$.
In view of Lemma~\ref{lem:support}, for all colours~$i$ dominated by~$c$, and thus for all colours~$i$ since~$c$ is dominating,
\[\Supp(W_{\bs{e_c}}^{DT}) \subseteq \Supp(W_{\bs{e_i}}^{DT}),\] 
implying that $\mathcal B(z_0, \varepsilon)\subseteq \bigcap_{i\in\{1,\ldots,d\}}\Supp(W_{\bs{e_i}}^{DT})$.
By Lemma~\ref{lem:densite},
we get that, for all $t\in[0,1]$,
\begin{equation}\label{eq:stable}
\big(t^{\nicefrac{\lambda}{S}} + (1-t)^{\nicefrac{\lambda}{S}}\big) \mathcal B(z_0, \varepsilon)
\subseteq \bigcap_{i\in\{1,\ldots,d\}}\Supp(W_{\bs{e_i}}^{DT}).
\end{equation}
Let us first assume that $\lambda\in \mathbb C\setminus \mathbb R$.
Note that when $t$ is tending to zero, 
$t^{\nicefrac{\lambda}{S}} + (1-t)^{\nicefrac{\lambda}{S}} = 1 + t^{\nicefrac{\lambda}{S}} + \mathcal O(t)$,
where, for all complex functions $\eta(t)$, we denote $\eta(t) = \mathcal O(t)$, if $|\eta(t)| = \mathcal O(t)$.
The function $t\mapsto t^{\nicefrac{\lambda}{S}}$ maps any neighbourhood of $0$ to an infinite spiral around~$1$ in the complex plane:
the module of $t^{\nicefrac{\lambda}{S}}$ goes to zero when $t$ goes to zero, 
while its argument, which is equal to $\mathtt{Im \lambda}\log t/S$ goes to $-\infty$.
In particular, the trajectory of $t\mapsto 1+ t^{\nicefrac\lambda{S}}$ 
intersects infinitely many times the set $\mathcal B(z_0, \varepsilon/2) \cap \{|z|=1\} \setminus\{1\}$.
Therefore, since $t^{\nicefrac{\lambda}{S}} + (1-t)^{\nicefrac{\lambda}{S}} = 1 + t^{\nicefrac{\lambda}{S}} + \mathcal O(t)$, when $t$ goes to zero, we get that there exists $t_0>0$ such that
\[|t_0^{\nicefrac{\lambda}{S}} + (1-t_0)^{\nicefrac{\lambda}{S}}| = 1,
\mathtt{arg}(t_0^{\nicefrac{\lambda}{S}} + (1-t_0)^{\nicefrac{\lambda}{S}}) \notin 2\pi\mathbb Z\text{ and } |t_0^{\nicefrac{\lambda}{S}} + (1-t_0)^{\nicefrac{\lambda}{S}}-1|< \varepsilon/2.\]
Therefore, the application
\[z\mapsto \big(t_0^{\nicefrac{\lambda}{S}} + (1-t_0)^{\nicefrac{\lambda}{S}}\big)z\]
is a rotation of an angle sufficiently small so that $\mathcal B(z_0, \varepsilon)$ intersects its image.
We also know, by Lemma~\ref{lem:densite} that $\bigcap_{i\in\{1,\ldots,d\}}\Supp(W_{\bs{e_i}}^{DT})$ is stable by this rotation.
Therefore, if we iterate this application, we get that there exists $\tilde \varepsilon>0$
$\{z\in\mathbb C \colon |z_0|- \tilde\varepsilon <|z|<|z_0|+\tilde\varepsilon\}$
is contained in $\bigcap_{i\in\{1,\ldots,d\}}\Supp(W_{\bs{e_i}}^{DT})$ (one can for example take $\tilde\varepsilon = \varepsilon/2$).
In particular, $-z_0 \in \bigcap_{i\in\{1,\ldots,d\}}\Supp(W_{\bs{e_i}}^{DT})$.
Using Lemma~\ref{lem:support}, we get that
\[\Big(\frac12\Big)^{\nicefrac{\lambda}{S}}\left(\mathcal B(z_0, \varepsilon) - z_0\right)
= \mathcal B \big(0,\left(\nicefrac12\right)^{\nicefrac{\RRe\lambda}{S}}\varepsilon\big)
\subseteq \bigcap_{i\in\{1,\ldots,d\}}\Supp(W_{\bs{e_i}}^{DT}).\]
Finally, recall that, in view of Lemma~\ref{lem:densite}, $\bigcap_{i\in\{1,\ldots,d\}}\Supp(W_{\bs{e_i}}^{DT})$ is stable by
\[z \mapsto (\nicefrac12)^{\nicefrac{\lambda}{S}-1} z,\]
which sends any ball $\mathcal B(0, r)$ to the 
larger ball $\mathcal B\big(0,(\nicefrac12)^{\nicefrac{\RRe\lambda}{S}-1}r\big)$ 
(recall that $\RRe\lambda < S$).
Iterating this application thus sends $\mathcal B \left(0,(\nicefrac12)^{\nicefrac{\RRe\lambda}{S}}\varepsilon\right)$ to $\mathbb C$, 
implying that $\bigcap_{i\in\{1,\ldots,d\}}\Supp(W_{\bs{e_i}}^{DT}) = \mathbb C$.

It remains to treat the case $\lambda\in \mathbb R$. 
From Proposition~\ref{thm:dcoul_dislocDT}, we have, for all $c\in\{1, \ldots, d\}$,
\[\sum_{i=1}^d(a_{c,i}+\delta_{c,i})\,\frac{\mathbb E W_{\bs e_i}^{DT}}{\theta_i+\lambda}=
\lambda\,\frac{\mathbb E W_{\bs e_c}^{DT}}{\theta_c+\lambda}.\]
First assume that, for all $1\leq i\leq d$, $\mathbb E W_{\bs e_i}^{DT} \geq 0$.
Then, if we take $c$ such that 
\[\frac{\mathbb E W_{\bs e_c}^{DT}}{\theta_c+\lambda} = \min_{1\leq i\leq d} \frac{\mathbb E W_{\bs e_i}^{DT}}{\theta_i+\lambda},\]
we get (see Equation~\eqref{eq:theta_i})
\[(S+1)\,\frac{\mathbb E W_{\bs e_c}^{DT}}{\theta_c+\lambda}
\leq \sum_{i=1}^d(a_{c,i}+\delta_{c,i})\,\frac{\mathbb E W_{\bs e_i}^{DT}}{\theta_i+\lambda}=
\lambda\,\frac{\mathbb E W_{\bs e_c}^{DT}}{\theta_c+\lambda},\]
implying that $\mathbb E W_{\bs e_c}^{DT} = 0$ since $\lambda < S$ by assumption.
The same reasoning holds if we assume that $\mathbb E W_{\bs e_i}^{DT} \leq 0$ for all $1\leq i\leq d$.
Therefore, either there exists $c$ such that $\mathbb E W_{\bs e_c}^{DT} = 0$,
or there exist $i$ and $j$ such that $\mathbb E W_{\bs e_i}^{DT} < 0 < \mathbb E W_{\bs e_j}^{DT}$.

First assume that exists~$c$ such that $\mathbb E W_{\bs e_c}^{DT} = 0$, 
and recall that, in view of Lemma~\ref{lem:esperance},
there exists a colour~$i$ such that $\mathbb E W_{\bs e_i}^{DT} \neq 0$. 
Therefore, there exists $z_i\neq 0$ in the support of $W_{\bs e_i}^{DT}$.
Since we have assumed that the urn is irreducible (assumption~$\mathtt{(I)}$), 
and thus, in particular that~$i$ is a dominating colour, 
we get, by Lemma~\ref{lem:support}, that $z_i\in \bigcap_{j=1}^d \Supp W_{\bs e_j}^{DT}$.
Similarly, the support of $W_{\bs e_c}^{DT}$, and thus $\bigcap_{j=1}^d \Supp W_{\bs e_j}^{DT}$, 
contain either $0$ or some $z_c$ chosen such that $z_i$ and $z_c$ have opposite signs. 
In the first case, since $W^{DT}_{\bs e_j}$ admits a continuous density for all $1\leq j\leq d$,
there exists $\varepsilon>0$ such that $[-\varepsilon, \varepsilon]\in \bigcap_{j=1}^d \Supp W_{\bs e_j}^{DT}$.
In latter case, applying Lemma~\ref{lem:support} 
implies that the segment $[z_i, z_c]$ (or $[z_c, z_i]$, depending on the signs) in
contained in $\bigcap_{j=1}^d \Supp W_{\bs e_j}^{DT}$. 
This segment contains~$0$, and thus, 
in both cases, there exists $\varepsilon>0$ such that
\[[-\varepsilon, \varepsilon]\subseteq \bigcap_{j=1}^d \Supp W^{DT}_{\bs e_j}.\]
Finally, by Lemma~\ref{lem:support}, $\bigcap_{j=1}^d \Supp W^{DT}_{\bs e_j}$ 
is stable by the application $z\mapsto (\nicefrac12)^{\nicefrac{\lambda}{S}-1}z$, 
and iterating this application sends $[-\varepsilon, \varepsilon]$ to~$\mathbb R$, which gives
\[\bigcap_{j=1}^d \Supp W^{DT}_{\bs e_j} = \mathbb R,\]
as desired.

It now remains to treat the case when there exist~$i$ and~$j$ 
such that $\mathbb E W_{\bs e_i}^{DT} < 0 < \mathbb E W_{\bs e_j}^{DT}$.
This is very similar to the previous case: there exists $z_i<0$ in $\Supp W^{DT}_{\bs e_i}$ 
and $z_j>0$ in $\Supp W^{DT}_{\bs e_j}$. By irreducibility, and applying Lemma~\ref{lem:support},
$[z_i, z_j]\subseteq \bigcap_{j=1}^d \Supp W^{DT}_{\bs e_j}$. This segment contains~$0$ in its interior 
and is therefore sent to~$\mathbb R$ by iterating the application $z\mapsto (\nicefrac12)^{\nicefrac{\lambda}{S}-1}z$,
under which $\bigcap_{j=1}^d \Supp W^{DT}_{\bs e_j}$ is stable.
\end{proof}

\begin{proof}[Proof of Theorem~\ref{th:density}]
Connection~\eqref{eq:mcCT} allows to transfer results stated in Proposition~\ref{prop:density} to $W_{\bs{e_c}}^{CT}$ for all $c\in\{1, \ldots, d\}$, since $\xi$ and $W_{\bs{e_c}}^{DT}$ both admit a density and are independent. Finally, Propositions~\ref{th:dcoul_decompDT} and~\ref{th:dcoul_decompCT} allow to generalise to any initial composition~$\bs\alpha$.
\end{proof}

\begin{remark}
For two-colour P{\'o}lya urns, it is proven under $\mathtt{(B)}$, $\mathtt{(I)}$ and $\mathtt{(T_{-1})}$ by \cite{CPS11} that
the density of $W^{CT}_{\bs \alpha}$ explodes in zero for all initial composition $\bs \alpha$. Is it possible to prove the same result for $d$-colour urns? We believe that the proof is similar to the one developed in~\cite{CPS11} if $\lambda$ is real, and under $\mathtt{(T)}$ instead of $\mathtt{(T_{-1})}$, but what can be said about $W^{CT}$ in zero when $\lambda$ is complex?
\end{remark}

\section{Conclusion and open problems}
This paper contains the first results concerning the $W$ random variables exhibited by the asymptotic behaviour of a multi--colour P{\'o}lya urn process projected along large Jordan blocks: we proved that $W$ is moments-determined and admits a density both in discrete and in continuous time.

However, many questions still remain about the random variables $W$: some of them are solved in the two--colour case (by methods that do not seem easy to generalise to multi--colour urns) but other are not proved even in the two--colour case. Can we find a lower bound for the moments of the $W$'s (\cite{CPS11} for two--colour urns)? What is the exact order of their moments? Is the Laplace transform of $W^{CT}$ convergent (this is false in the two--colour case~\cite{CPS11})? Is the Fourier transform of $W^{CT}$ integrable (also false in the two--colour case~\cite{CPS11})? Is the density of $W^{CT}$ continuous (also false in two--colour~\cite{CPS11})? Do the $W$'s have heavy tails?

\section*{Acknowledgements} The author is grateful to Brigitte Chauvin, Svante Janson and Nicolas Pouyanne for fruitful discussions about these results. She would also like to warmly thank the anonymous referee for the careful proof-reading of the preliminary versions of this paper, and the EPSRC for support through the grant EP/K016075/1.

\bibliographystyle{alpha}
\bibliography{urnes}
\end{document}